\documentclass[leqno]{article}
\usepackage[hmargin=30mm,vmargin=30mm]{geometry}
\usepackage{amsmath, amsthm, amssymb, mathabx, verbatim, setspace, enumerate,mathtools, caption}
\usepackage[mathscr]{euscript}
\usepackage[all,cmtip]{xy}
\usepackage{etoolbox}
\usepackage{hyperref}
\usepackage[OT2,T1]{fontenc}

\newcommand{\CC}{\mathbb{C}}
\newcommand{\ds}{\displaystyle}
\newcommand{\ra}{\rightarrow}
\newcommand{\ZZ}{\mathbb{Z}}

\newcommand{\QQ}{\mathbb{Q}}

\newcommand{\OO}{\mathcal{O}}

\newcommand{\HH}{\mathcal{H}}

\newcommand{\C}{\mathcal{C}}
\newcommand{\A}{\mathcal{A}}

\newcommand{\SSS}{\mathcal{S}}

\newcommand{\AAA}{\mathbb{A}}
\newcommand{\pp}{\mathfrak{p}}

\newcommand{\eps}{\varepsilon}

\newtheorem*{Thm*}{Theorem}
\newtheorem*{cor*}{Corollary 1.2'}

\newtheorem{Thm}{Theorem}[section]
\newtheorem{Prop}[Thm]{Proposition}
\newtheorem{Lem}[Thm]{Lemma}
\newtheorem{cor}[Thm]{Corollary}

\newtheorem{rmk}[Thm]{Remark}

\DeclareMathOperator{\GL}{GL}

\DeclareMathOperator{\sgn}{sgn}

\DeclareMathOperator{\SL}{SL}
\DeclareMathOperator{\Gal}{Gal}
\DeclareMathOperator{\Nm}{Nm}

\DeclareMathOperator{\Char}{char}

\DeclareMathOperator{\vol}{vol}

\DeclareSymbolFont{cyrletters}{OT2}{wncyr}{m}{n}
\DeclareMathSymbol{\Sha}{\mathalpha}{cyrletters}{"58}

\DeclareMathOperator{\Cl}{Cl}
\DeclareMathOperator{\Tr}{Tr}

\begin{document}
\title{$L$-series values of sextic twists for elliptic curves over $\QQ[\sqrt{-3}]$}  
\author{Eugenia Rosu}
\date{}
\maketitle

\begin{abstract} We prove a new formula for the central value of the $L$-function $L(E_{D, \alpha}, 1)$ corresponding to the family of sextic twists over $\QQ[\sqrt{-3}]$ of elliptic curves $E_{D, \alpha}: y^2=x^3+16D^2\alpha^3$ for $D$ an integer and $\alpha \in \QQ[\sqrt{-3}]$. The formula generalizes the result of cubic twists over $\QQ$ of Rodriguez-Villegas and Zagier for a prime $D \equiv 1 (9)$ and of Rosu for general $D$. For $\alpha$ prime and all integers $D$, we also show that the expected value from the Birch and Swinnerton-Dyer conjecture of the order of the Tate-Shafarevich group is an integer square in certain cases, and an integer square up to a factor $2^{2a}3^{2b}$ in general.
\end{abstract}


\section{Introduction}

  In this paper we will consider the family of elliptic curves:
\begin{equation}\label{main}
E_{D, \alpha}: y^2=x^3+16D^2\alpha^3,\ \  D \in \ZZ, \ \alpha \in \ZZ[\omega],
\end{equation}
where $\omega=\frac{-1+\sqrt{-3}}{2}$. The elliptic curves $E_{D, \alpha}$ are quadratic twists over $K=\QQ[\sqrt{-3}]$ of the elliptic curve $y^2=x^3-432D^2$, which is a Weierstrass equation for the famous family of elliptic curves 
\begin{equation}\label{cubes}
E^{(D)}: x^3+y^3=D.
\end{equation}
In particular, $E_{D, \alpha}$ are sextic twists over $K$ of the familiar elliptic curve $E_1: y^2=x^3+1$. The equation \eqref{cubes} was extensively studied in the literature over $\QQ$, with the goal to answer the question for which integers we can write $D$ as the sum of two rational cubes.

 However, very little is known about the twists of $E^{(D)}$ over $K$. The goal of this paper is to study the family of twists $E_{D, \alpha}$ and, in particular, the central values $L({E_{D, \alpha}}/K, 1)$ of their $L$-functions. From the {\it Birch and Swinnerton-Dyer(BSD)}  conjecture, the vanishing of $L({E_{D, \alpha}}/K, 1)$ is equivalent to having rational solutions for $E_{D, \alpha}(K)$. Without assuming BSD, when $L({E_{D, \alpha}}/K, 1)\neq 0$ we have no rational solutions $E_{D, \alpha}(K)$ from the work of Coates and Wiles \cite{CW}. 

We define the invariant
\begin{equation}
S_{D\alpha}=\frac{1}{c_{E_{D, \alpha}}\Omega_{D\alpha}\overline{\Omega}_{D\alpha}}L({E_{D, \alpha}}/K, 1),
\end{equation}
where $c_{E_{D, \alpha}}=\prod\limits_{v|6D\alpha}c_v$ is the product of Tamagawa numbers $c_v$ and $\Omega_{D\alpha}\in \CC^{\times}$ is a period of $E_{D, \alpha}$, more precisely $\Omega_{D\alpha}\overline{\Omega}_{D\alpha}=\frac{1}{m^{1/2}}(\frac{\Gamma(1/3)^3}{4\pi D^{1/3}})^2$. We have defined $S_{D\alpha}$ such that when $\Nm(\alpha)>1$, if $L(E_{D, \alpha}/K, 1)\neq 0$, from the BSD conjecture we have:
\[
S_{D\alpha}=\#\Sha_{E_{D, \alpha}/K},
\] 
where $\Sha_{E_{D, \alpha}/K}$ is the order of the Tate-Shafarevich group of $E_{D, \alpha}$ over $K$. We note that we have used here the fact that the torsion subgroup $E_{D, \alpha}(K)_{tor}$ is trivial for $\Nm\alpha>1$.

When $L({E_{D, \alpha}}/K, 1)\neq 0$, the order of $\Sha_{E_{D, \alpha}}$ over $K$ is known to be finite from the work of Rubin \cite{Ru}, and its order is a square as proved by Cassels \cite{C} via the Cassels-Tate pairing. We will show in the current paper that, for certain cases, $S_{D\alpha}$ is indeed an integer square, and in general $16 S_{D\alpha}$ is an integer square up to an even power of $3$. 

By computing $S_{D\alpha}$, we can check if we have solutions in $\eqref{main}$ and, if $S_{D\alpha}\neq 0$, we get the value of the analytic rank $S_{D\alpha}$ of the Tate-Shafarevich group. To summarize:

  \begin{enumerate}[(i)]
   	\item $S_{D\alpha}\neq 0$ $\Longrightarrow$ no solutions in \eqref{main}
	\item $S_{D\alpha}\neq 0$ $\xRightarrow{BSD}$ $S_{D\alpha}=\#\Sha$ integer square
	\item  $S_{D\alpha}=0$ $\xRightarrow{BSD}$ have solutions in \eqref{main}.
  \end{enumerate}
  
  From the work of Rubin \cite{Ru}, it is known that $v_p(S_{D\alpha})=v_p(\Sha_{E_{D, \alpha}/K})$ for all primes $p\neq 2, 3$, where $v_p$ is the valuation at $p$. We show that further $v_p(S_{D\alpha})\equiv v_p(\Sha_{E_{D, \alpha}/K})(2)$ for $p=2, 3$.

We also note that the $L$-function $L(E_{D, \alpha}/K, s)$ has a functional equation from $s$ to $2-s$ that is expected to have constant global root number $w=1$ (see \cite{BK}), thus we do not have a priori expectations for any of the central values $L(E_{D, \alpha}/K, 1)$ to vanish.


\bigskip

In the current paper, we will compute several formulas for $S_{D\alpha}$. In previous work, Rodriguez-Villegas and Zagier computed in \cite{RV-Z} the special values $L(E^{(p)}, 1)$ of $L$-functions over $\QQ$ for primes $p\equiv 1(9)$, and Rosu extended their work to all integers $D$ in \cite{Ro}. We generalize their work to compute $L(E_{D, \alpha}/K, 1)$ in Theorem \ref{thm1} and use novel methods to show that these values are squares in Theorem \ref{thm2}. We note that the case of quadratic twists over $\QQ$ was widely studied using such powerful tools as Waldspurger's theorem (\cite{W}). However these tools cannot be applied in our case as we are considering quadratic twists over the imaginary quadratic field $K$ as well as cubic twists by $D$. Instead we use the properties of the elliptic curves $E_{D, \alpha}$ that have complex multiplication (CM) by $\OO_K$, the ring of integers of the number field $K$.

\bigskip

We will state now the results of the paper. Without loss of generality, we choose $\alpha$ square-free with $\Nm\alpha>1$ and $D$ cube-free. Moreover, we note that the elliptic curve $E_{D, \alpha}$ is invariant when multiplying $\alpha$ by a cubic root of unity, thus we will fix throughout the paper the representatives 
\[
\alpha \equiv \pm 1, \pm \sqrt{-3} (4).
\]
We denote $m=\Nm\alpha$.


\begin{Thm}\label{thm3} For $\alpha$ prime and any integer $D$ such that $(D, 6\alpha)=(\alpha, 6)=1$,  $S_{D\alpha}$ is an integer square up to a factor $2^{2a}\cdot 3^{2b}$, $0\leq a\leq 2$. More precisely:
	\begin{itemize}
		\item $v_p(S_{D\alpha})$ is even for all $p$, 
		\item  $v_p(S_{D\alpha})\geq 0$, for all $p\neq 2, 3$ and 
		\item $v_3(c_{E_{D, \alpha}}S_{D\alpha})\geq -1$,  $v_2(c_{E_{D, \alpha}}S_{D\alpha})\geq 2e$, where $e=1$ for $m\equiv1(4)$ and $e=0$ for $m\equiv 3(4)$.
	\end{itemize}	 
\end{Thm}

Moreover, as a particular case, $S_{D\alpha}$ is an integer square for $c_{E_{D, \alpha}}=1, 4$ when $m\equiv 1(4)$, or $c_{E_{D, \alpha}}=1$ when $m\equiv 3(4)$. Taking $\left[\frac{\alpha}{\cdot}\right]$ to be the quadratic character over $K$ and $\chi_{2D^2}$ the cubic character (see Section \ref{hecke} for the definitions), then we have:

\begin{Thm}\label{cor1} $S_{D\alpha}$ is an integer square for $\alpha$ prime with $m\equiv 1(4)$ such that

\begin{itemize}
	
	\item $D\equiv \pm 1(9)$, $\alpha\equiv -1(3)$ and $\left[\frac{\alpha}{\pp}\right]=-1$ for all prime ideals $\pp|D$
	
	\item $D\equiv \pm 4(9)$, $\alpha\equiv 1(3)$ and $\left[\frac{\alpha}{\pp}\right]=-1$ for all prime ideals $\pp|D$
	
	\item $D\equiv \pm 2(9)$ and $\chi_{2D^2}(\alpha)\neq 1$
\end{itemize}
	
	 For $m\equiv 3(4)$, $S_{D\alpha}$ is an integer square for $D\equiv \pm 1, \pm 4$ such that $\chi_{2D^2}(\alpha)\neq 1$ and  satisfying the conditions above.
	 
\end{Thm}


	


Moreover, when $D=1$, for the elliptic curve $E_{\alpha}: y^2=x^3+16\alpha^3$ defined over $K$, we have:
\begin{cor}\label{cor2} For $\alpha\equiv -1(3)$ prime, $S_{\alpha}$ is an integer square. 
\end{cor}

 Theorem \ref{thm3} follows from the explicit formula of $S_{D\alpha}$ that we compute in Theorem \ref{thm2}:

\begin{Thm}\label{thm2} Let $\alpha$ prime and  $D$ any integer such that $(\alpha, 3)=1$, $(D, 6\alpha)=1$. Then for $c_{E_{D, \alpha}}=4^a\cdot 3^b$, $0\leq b\leq 2$, the Tamagawa number of $E_{D, \alpha}$, we have
\[ 
S_{D\alpha}=\frac{1}{4^{a}(-3)^{b+1}} Z^2,  
\]
where $Z$ is an integer if $b$ is odd, and $Z/\sqrt{-3}$ is an integer if $b$ is even.

Here 
\[
Z=u\Tr_{H_{3D^*}/K}\frac{\Theta_M(D^*\tau_0/m^*)}{\Theta_K(\tau_0/m)}D^{1/3} \alpha^{1/2},
\] 
\begin{itemize} 
	
	\item $\tau_0=\frac{b_0+\sqrt{-3}}{2}$ CM point with $b_0\equiv \sqrt{-3} (\alpha)$, $b_0\equiv 1(2)$

	\item $M=\QQ[\sqrt{-3m}]$ of discriminant $-3m^*$ and class number $h_M$,  $H_{3m^*D}$ is the ring class field over $K$ corresponding to the order $\OO_{3Dm^*}=\ZZ+3Dm^*\OO_K$
	\item $\Theta_M(z)=h_M+2\sum\limits_{N\geq 1}\sum\limits_{d|N}\left(\frac{d}{3m^*}\right)e^{2\pi i Nz}$ and $\Theta_K(z)=\sum\limits_{m, n\in \ZZ} e^{2\pi i (m^2+n^2-mn)z}$ theta functions of weight $1$

	\item $D^*=\begin{cases} D & \alpha\equiv 1(4) \\ 4D & \alpha\equiv -1, \pm\sqrt{3}(4)\end{cases}$, $u=2^e\omega^l$, for $\omega^l$ a cubic root of unity, and $e=\begin{cases}1 & m\equiv 1(4)\\0 & m\equiv 3(4)\end{cases}$

	\end{itemize}
		
\end{Thm}

\bigskip

 We get immediately Theorem \ref{thm3}, as well as Theorem \ref{cor1} and Corollary \ref{cor2}. We also note the surprisingly simple formulas for $\alpha$ prime:
 
\begin{cor} For $\alpha$ prime such that $\alpha\equiv 1(4)$, we have \[
S_{\alpha}=\left(\frac{2}{\sqrt{3}}\frac{\Theta_M(\tau_0/m)}{\Theta_k(\tau_0/m)}\alpha^{1/2}\right)^2.
\] 
For $\alpha$ prime such that $\alpha\equiv -1, \pm \sqrt{-3}(4)$, we have $S_{\alpha}=\left(\frac{2}{\sqrt{3}}
\Tr_{H_{12}/K}\frac{\Theta_M(\tau_0/m)}{\Theta_k(\tau_0/m)}\alpha^{1/2}\right)^2$.
\end{cor}

As mentioned in Corollary \ref{cor2}, for $\alpha\equiv -1(3)$, $S_{\alpha}$ is an integer square. Similarly, in the case  $\alpha\equiv 1(3)$, $9S_{\alpha}$ an integer square.
 \bigskip
 
  Theorem \ref{thm2} is based on the more general result Theorem \ref{thm1}, proved for all $\alpha$ not necessarily prime. From CM theory we can find a Hecke character $\chi$ defined over $\AAA^{\times}_K$, the ideles of $K$, such that:
\[
L({E_{D, \alpha}}/K, s)=L(s, \chi)L(s, \overline{\chi}).
\]
By computing the value of each $L(1, \chi)$, we get:

\begin{Thm}\label{thm1} For $D\in \ZZ$ and $\alpha\in \OO_K$ such that $(D, 6\alpha)=1$, $(\alpha, 3)=1$, we have:
\[
S_{D\alpha}=\frac{1}{c_{E_{D, \alpha}}}\left|c\Tr_{H_{3D'm^*}/K} \frac{\Theta_M(D'\omega)}{\Theta_K(\omega)} \alpha^{1/2}D^{1/3}\right|^2,
\]

\noindent and $ S_{D\alpha}/|c|^2\in \ZZ$. Here $c=\frac{2L_{\overline{\alpha}}(1, \chi)|}{\sqrt{3m^*}}$, where $L_{\overline{\alpha}}(1, \chi)=\prod_{\pp|(\alpha)}(1-\chi(\pp)\Nm\pp^{-1})^{-1}$, and $D'=\begin{cases} D& \alpha \equiv 1, \pm\sqrt{-3} (4) \\ 4D & \alpha\equiv -1(4)\end{cases}$.

\end{Thm}

\bigskip

The paper is structured as follows. In Section \ref{hecke} we present the background on the Hecke characters, in particular the properties of the quadratic character $\left[\frac{\alpha}{\cdot} \right]$. In the rest of Section \ref{background} we cover the properties of the weight $1$ classical Eisenstein series and prove several properties of the theta function $\Theta_M(z)$, including a Siegel-Weil type result and an inverse transformation. 

The goal of Section \ref{L} and Section \ref{Galois} is to prove Theorem \ref{thm1}. The proof is similar to that of Theorem 3 of \cite{Ro} and consists of computing a formula for $L(s, \chi)$ by using Tate's thesis, which we cover in section \ref{L}. We obtain a finite linear combination of Eisenstein series and characters and by further applying the Siegel-Weil type result we get a linear combination of ratios of theta functions evaluated at CM points. These CM points correspond to the ideal classes of the ring class group $\Cl(\OO_{3D'm^*})$. We finish up the proof in Section \ref{Galois} by applying Shimura reciprocity law to show that all terms are conjugate to each other over $H_{3D'm^*}$.

Section \ref{X} is the heart of the paper. The goal is to show Theorem \ref{thm2} by rewriting the term $X_{D\alpha}=\Tr_{H_{3D'm^*}/K} \frac{\Theta_M(D'\omega)}{\Theta_K(\omega)}D^{1/3} \alpha^{1/2}$ of  Theorem \ref{thm1}. We prove in Section \ref{XT} that $L_{\overline{\alpha}}(1, \chi)X_{D\alpha}/\alpha$ equals $Z^{\circ}=\Tr_{H_{3D^*}/K}\frac{\Theta_M(D^*\tau_0/m^*)}{\Theta_K(\tau_0/m)}D^{1/3} \alpha^{1/2}$ up to a cubic root of unity. In section \ref{XYZ}, we relate $Z^{\circ}$ to its complex conjugate and show that $Z^{\circ}\in \OO_K$ is either real or purely imaginary. These two steps require relating $X_{D\alpha}$ and $Z^{\circ}$ to the traces of various modular functions using properties of the theta function $\Theta_M$ and Shimura reciprocity law applied to several Galois conjugates. The method should work for general $m$ and $D$, however this is very technical and quite involved. In order to simplify the computations, we treated only the case of $\alpha$ prime in Theorem \ref{thm2}.

In section \ref{comp_shimura} we cover several technical results proved using Shimura reciprocity law that are used in Section \ref{X}. Finally, in the Appendix we present the explicit computations involving the Galois conjugates of modular functions used in the proofs of Theorem \ref{thm2}. These are proved using Shimura's reciprocity law and properties of the theta function $\Theta_M$ proved in Section \ref{theta_M}. We also compute in this section the explicit values of the cubic and quadratic characters for various ideals.

\bigskip

We note that one can computationally check the values of $S_{D\alpha}$. These computations, while delicate to program, are much faster than the direct computations of the $L$-functions in Magma.


\bigskip
{\bf Acknowledgements.} The author would like to thank the Max Planck Institute for Mathematics for hospitality. The author was also partially supported by the grant CRC 326 - Gaus. We thank Don Zagier for interesting conversations.

\section{Background}\label{background}

Let $K=\QQ[\sqrt{-3}]$, $\OO_K=\ZZ[\omega]$ its ring of integers for $\omega=\frac{-1+\sqrt{-3}}{2}$ a fixed cubic root of unity, and $K_v$ the completion of $K$ at the place $v$. We denote $K_p=K\otimes_{\QQ}\QQ_p$ and $\OO_{K_p}=\OO_{K}\otimes_{\ZZ} \ZZ_p$ the semilocal ring of integers.

We also remark that we can write the primitive ideals $\A$ in $\OO_K$ as $\ZZ$-modules $\A=\left[a, \frac{-b+\sqrt{-3}}{2}\right]$, where $a=\Nm\A$ and $b$ is chosen (non-uniquely) such that $b^2\equiv -3 (4a)$. We will denote $\tau_{\A}=\frac{-b+\sqrt{-3}}{2a}$.

\subsection{Hecke characters}\label{hecke}

We recall that from CM theory we can a find a Hecke character $\chi:\AAA^{\times}_K \ra \CC^{\times}$ corresponding to the elliptic curve $E_{D, \alpha}$ such that:
\begin{equation}
L({E_{D, \alpha}}/K, s)=L(s, \chi)L(s, \overline{\chi}).
\end{equation}
We can explicitly write the Hecke character $\chi$ as a product of Hecke characters $\varphi, \chi_D, \varepsilon$ defined over $K$:
\begin{equation}\label{chi}
\chi=\varphi \chi_D\varepsilon.
\end{equation}

We define the characters $\varphi, \chi_D, \varepsilon$ in classical language. The character $\varphi$ has conductor $3$ and, for an ideal $\A$ prime to $3$, it is defined by 
\[
\varphi(\A)=k_{\A},
\]
for $k_{\A}$ the unique generator of $\A$ such that $k_{\A}\equiv 1(3)$. As $K$ is a PID and $\OO_K^*$ is generated by $-\omega$, the character is well-defined.

\subsubsection{The cubic character $\chi_D$}\label{cubic} The cubic character $\chi_D: I(3D) \ra \{1, \omega, \omega^2\}$ is defined for the space of ideals prime to $3D$ to be $\chi_D=\overline{\left(\frac{D}{\cdot}\right)}_3$, the complex conjugate of the cubic character defined in Ireland and Rosen \cite{IR}. More precisely, $\left(\frac{D}{\cdot}\right)_3$ is the unique cubic root of unity for which we have $\left(\frac{D}{\pp}\right)_3\equiv D^{\frac{\Nm\pp -1}{3}} \mod  \pp$.  Moreover, $\chi_D$ is well-defined on $\Cl(\OO_{3D})$,  the class group for the order $\OO_{3D}=\ZZ+3D\OO_K$.
We note the property:
\[
(D^{1/3})^{\sigma_{\A}^{-1}}= D^{1/3} \chi_D(\A),
\]
where $\sigma_{\A}$ is the Galois action corresponding to the ideal $\A$ via the Artin map.

If $\beta\in \OO_K$ such that $\beta \equiv \omega (D)$, $\beta\equiv 1 (3)$ we will write $\chi_D(\beta)=\chi_{\omega}(D)$ and this equals $\overline{\omega}^{\frac{D^2-1}{3}}$ for $D\equiv 2(3)$, and $\overline{\omega}^{2\frac{D-1}{3}}$ for $D\equiv 1(3)$, respectively. Thus explicitly we have:
\[
\chi_{\omega}(D)=\begin{cases}
1 & \text{ if } D\equiv \pm 1 (9) \\
\omega & \text{ if } D\equiv \pm 4 (9) \\
\omega^2 & \text{ if } D\equiv \pm 2 (9).\\
\end{cases}
\]


\subsubsection{The quadratic character $\eps$}\label{quad_char}	The quadratic character $\varepsilon: I(2\alpha) \ra \{\pm 1\}$ is defined on the space of ideals prime to $2\alpha$ and it  has conductor $\alpha$ for $\alpha\equiv 1 (4)$ and conductor $4\alpha$ for $\alpha\equiv -1, \pm\sqrt{-3}(4)$. We have $\eps=\left[\frac{\alpha}{\cdot}\right]$, where $\left[\frac{\alpha}{\cdot}\right]$ is defined as in  \cite{IR} by taking 
	\[
	\left[\frac{\alpha}{\pp}\right] \equiv \alpha^{\frac{\Nm\pp - 1}{2}} \mod \pp
	\]
for a prime ideal $\pp$ prime to $2\alpha$, and extending by multiplication. We also write $\left[\frac{\alpha}{\beta}\right]=\left[\frac{\alpha}{(\beta)}\right]$ and we note the property:
\[
(\alpha^{1/2})^{\sigma_{\A}^{-1}}= \alpha^{1/2} \eps(\A),
\]
where $\sigma_{\A}$ is the Galois action as above. 
We note $\left[\frac{\alpha}{\beta}\right]=\left[\frac{\overline{\alpha}}{\overline{\beta}}\right]$, thus $\left[\frac{\alpha}{\A}\right]=\left[\frac{\overline{\alpha}}{\overline{\A}}\right]$.

\bigskip

We will follow Lemmermeyer \cite{L2} in presenting several properties of the quadratic character $\left[\frac{\alpha}{\cdot}\right]$ for $\alpha\in \QQ[\sqrt{-3}]$, in particular the Eisenstein quadratic reciprocity law. We note that $\omega$ in the notation of \cite{L2} is $-\omega^2$ in our notation.

Clearly $\left[\frac{\alpha}{\pp}\right]=\left[\frac{\alpha'}{\pp}\right]$ if $\alpha \equiv \alpha' \mod \pp$, and, moreover, for $a\in \ZZ$ and $\beta\in \OO_K$, we have:
\begin{equation}\label{quad_norm}
\left[\frac{\beta}{a}\right]=\left(\frac{\Nm\beta}{a}\right), \ \ \ \left[\frac{a}{\beta}\right]=\left(\frac{a}{\Nm\beta}\right),
\end{equation}
where $\left(\frac{\cdot}{a}\right)$,  $\left(\frac{a}{\cdot}\right)$ are the usual quadratic characters over $\QQ$.

\bigskip

{\it Eisenstein's quadratic reciprocity law} for $\alpha, \beta$ of odd norm is given by:
\begin{equation}\label{rec_law0}
\left[\frac{\alpha}{\beta}\right]=
\left[\frac{\beta}{\alpha}\right][\alpha][\overline{\beta}][\alpha\overline{\beta}],
\end{equation}
for a symbol $[\gamma]$ that is uniquely determined by the value of $\gamma$ modulo $4$. For general $\gamma=a-b\omega^2$ of odd norm, we take as definition $[\gamma]=(-1)^{\frac{\Nm\gamma -1}{2} \frac{b'-1}{2}} \left(\frac{2}{\Nm\gamma}\right)^r,$ where $b=2^rb'$ for $b'$ odd (see Lemma 12.11 of \cite{L2}). However, for $\gamma=a-b\omega^2$ a primitive element of odd norm, the definition simplifies to
\[
[\gamma]=\left(\frac{b}{\gamma}\right).
\] 

We have the supplementary reciprocity laws:
\begin{enumerate}[(i)]
	\item $\left[\frac{-1}{\gamma}\right]=(-1)^{\Tr\frac{\gamma-1}{2}},$
	\item $\left[\frac{2}{\gamma}\right]=(-1)^{\Tr\frac{\gamma^2-1}{8}}.$

\end{enumerate}

We note that the usual quadratic reciprocity law holds for $\alpha\equiv 1(4)$. For the remaining cases $\alpha \equiv -1, \pm\sqrt{-3}$ the reciprocity law does not hold in general, but depends on $\beta \mod 4$. Noting that we can modify $\beta$ by a cubic root of unity to get $\beta\equiv \pm 1, \pm \sqrt{-3}(4)$, we compute explicitly in the following lemma:

\begin{Lem}\label{rec_law} For $\beta$ prime to $2\alpha$, we have:
\begin{enumerate}[(i)] 
	\item For $\alpha\equiv 1(4)$, we have the reciprocity law $\left[\frac{\alpha}{\beta}\right]=
\left[\frac{\beta}{\alpha}\right]$.

	\item For $\alpha\equiv -1(4)$, we have:
	$\left[\frac{\alpha}{\beta}\right]
\left[\frac{\beta}{\alpha}\right]=\begin{cases}
		1 & \text{for }\beta \equiv \pm 1 (4) \\ 
		\left[\frac{\beta/\sqrt{-3}}{4}\right] &\text{for } \beta \equiv \pm \sqrt{-3} (4)
\end{cases}$

	\item For $\alpha\equiv \pm\sqrt{-3}(4)$, we have $\left[\frac{\alpha}{\beta}\right]\left[\frac{\beta}{\alpha}\right]=\begin{cases} 
		\left[\frac{\beta}{4}\right] & \text{for }\beta \equiv \pm 1 (4) \\ 
		-\left[\frac{\alpha/\sqrt{-3}}{4}\right]\left[\frac{\beta/\sqrt{-3}}{4}\right] & \text{for }\beta \equiv \pm \sqrt{-3} (4)
		 \end{cases}.$
\end{enumerate}
\end{Lem}

{\bf Proof:} We compute $[1]=[-1]=[-\sqrt{-3}]=1$ and $[\sqrt{-3}]=-1$ and apply \eqref{rec_law0}.

\bigskip

We get immediately from Lemma \ref{rec_law} and \eqref{quad_norm}: 

\begin{cor}\label{place_2} For an ideal $\A$ with generator $k_{\A}$ such that $\begin{cases} k_{\A}\equiv n (m) & \text{for }\alpha\equiv 1 (4) \\ k_{\A}\equiv n (4m) & \text{for } \alpha\equiv -1, \pm \sqrt{-3} (4)\end{cases}$, with $n\in \ZZ$, we have
\[
\eps(\A)=\left(\frac{n}{m^*}\right).
\]
\end{cor}

Corollary \ref{place_2} will be used in Section \ref{L} when choosing Schwartz-Bruhat functions $\Phi_p$ such that $\chi_p \Phi_{p}$ to be invariant on $(\ZZ+3D'm^*\OO_{K_p})$ for $p|3D'm^*$.

We further present the following lemma that will be used in Section \ref{XT}:

\begin{Lem}\label{eps_1} For $\alpha \equiv 1$ primitive, we have $\eps(\overline{\alpha})=\eps(\sqrt{-3})$ and for $\alpha\equiv -1, \pm\sqrt{-3}(4)$ primitive, we have $\eps(\overline{\alpha})=-\eps(\sqrt{-3})$.
\end{Lem}
{\bf Proof:} By definition, $[\alpha]=\left[\frac{(\alpha-\overline{\alpha})/(\sqrt{-3})}{\alpha}\right]$, which equals $[\alpha]=\left(\frac{-1}{m}\right)\left[\frac{\overline{\alpha}}{\alpha}\right]\left[\frac{\sqrt{-3}}{\alpha}\right]$. By the reciprocity law \eqref{rec_law0}, we have $\left[\frac{\sqrt{-3}}{\alpha}\right]=\left[\frac{\alpha}{\sqrt{-3}}\right][\alpha][-\sqrt{-3}][-\sqrt{-3}\alpha]$, thus $[\alpha]=\left(\frac{-1}{m}\right)\left[\frac{\overline{\alpha}}{\alpha}\right]\left[\frac{\alpha}{\sqrt{-3}}\right][\alpha][-\sqrt{-3}\alpha]$, which implies $\left[\frac{\overline{\alpha}}{\alpha}\right]=\left(\frac{-1}{m}\right)\left[\frac{\alpha}{\sqrt{-3}}\right][-\sqrt{-3}\alpha]$. But $[-\sqrt{-3}\alpha]=1$ for $\alpha \equiv 1, \pm\sqrt{-3} (4)$ and $[-\sqrt{-3}\alpha]=-1$ for $\alpha \equiv -1(4)$. Thus $\eps(\sqrt{-3})=\eps(\overline{\alpha})$ for $\alpha\equiv 1(4)$ and $\eps(\sqrt{-3})=-\eps(\overline{\alpha})$ for $\alpha\equiv \pm \sqrt{-3}(4), -1$.

\subsection{Adelic Hecke characters}

Each Hecke character $\chi$ can be written adelically as $\chi=\prod_{v} \chi_v$, where $\chi_v:K_v^{\times}\ra \CC$, such that at the unramified places we have $\chi_v(\OO_{K_v}^{\times})=1$ and $\chi_v(\varpi_v)=\chi(\pp_v)$, where $\varpi_v$ is the uniformizer corresponding to the prime ideal $\pp_v$. 

In our case, we have:
\begin{itemize}
	\item $\varphi_{\infty}(z)=z^{-1}$, $\varphi_p(p)=-p$ for $p\equiv 2(3)$, $\varphi_v(\varpi_v)=\varpi_v$, where $\varpi_v\equiv 1(3)$ uniformizer for $\pp_v$ prime such that $\pp_v\overline{\pp_v}=p$ prime in $\QQ$
	
	\item $\chi_{D, \infty}(z)=1$, $\chi_{D, p}(p)=1$ for $p\equiv 2(3)$ and $\chi_{D, v}(p)\chi_{D, \overline{v}}(p)=1$ for $\pp_v\pp_{\overline{v}}=p$ prime in $\QQ$
	
	\item $\eps_{\infty}(z)=1$, $\eps_{v}(\omega_v)=\left(\frac{\alpha}{\pp_v}\right)$

\end{itemize}

\subsection{Eisenstein series}\label{eisenstein}
In the current section we will discuss the properties of a family of weight $1$ Eisenstein  series $E_N(s, z)$ and relate their central values $E_N(0, z)$ to theta functions $\Theta_N(z)$. The Eisenstein series $E_N(0, z)$ will appear in Section \ref{L} in the computation of the $L$-function $L(1, \chi)$ for $N=3m^*$ and we will discuss further properties of the theta function $\Theta_M$ corresponding to $E_{3m^*}$ in the following section. 

For a positive integer $N$, we define the following classical Eisenstein series:
\[
E_s(z; a_1, a_2, N)=\sum_{\substack{c\equiv a_1 (\text{mod~} N) \\ d\equiv a_2 (\text{mod~} N) }} \frac{1}{(cz+d)|cz+d|^s}, 
\] 
for $z\in \HH$ and $s\in \CC$. They were extensively studied by Hecke in \cite{H}. While the series does not converge absolutely for $s=0$, we can still compute their Fourier expansion using Hecke's trick. More precisely, Hecke computed: 
\begin{equation}\label{fourier_eis}
E_0(z; a_1, a_2, N)= C(a_1, a_2, N)-\frac{2\pi i}{N} \sum\limits_{\substack{cn>0, \\ c\equiv a_1 (\text{mod } N)}} (\sgn (n)) e^{2\pi i na_2/N} e^{2\pi i nc z/N},
\end{equation}
\noindent where the constant part is $C(a_1, a_2, N)=\delta(a_1/N)\sum\limits_{d \equiv a_2(\text{mod } N)} \frac{\sgn(d)}{|d|^{s+1}}|_{s=0}-\frac{\pi i }{N} \sum\limits_{c \equiv a_1 (\text{mod } N)} \frac{\sgn(c)}{|c|^s}|_{s=0}$, for $\delta(a_1/N)=\begin{cases} 1, & \text{ if } N|a_1; \\ 0, & \text{if } N\nmid a_1.   \end{cases}$

We are interested in the Eisenstein series
\begin{equation}
E_N(s, z)=\sum_{a_2=0}^{N-1} \left(\frac{a_2}{N}\right)E_{s}(z; 0, a_2, N).
\end{equation}
 We will relate its central value $E_N(0, z)$ to the theta function:
\begin{equation}
\Theta_N(z)=h_N+2\sum\limits_{n\geq 1} (\sum\limits_{d|n} \left(\frac{d}{N}\right))e^{2\pi i n z},
\end{equation}
where $-N$ is the discriminant of the number field $L=\QQ[\sqrt{-N}]$ and $h_N$ is the class number of $\OO_L$.

More precisely, comparing the Fourier expansions, we compute below a variation of the Siegel-Weil theorem: 

\begin{Lem}\label{Theta_Eis}For $N\not\eq 3$, we have $E_{N}(0, z)=
\frac{2 \pi}{\sqrt{N}}\Theta_{N}(z)$.

\end{Lem}
{\bf Proof:}  As $\ds E_{N}(0, z)=\sum_{a_2=0}^{N-1} \left(\frac{a_2}{N}\right) E_0(z; 0, a_2, N)$, we use the Fourier expansion \eqref{fourier_eis} at $s=0$ for each $E(z; 0, a_2, N)$. We compute first the constant term: 
\[
\ds \sum_{a_2=0}^{N-1}\left(\frac{a_2}{N}\right) C(0, a_2, N)
=
2\sum_{d=1}^{\infty} \frac{\left(\frac{d}{N}\right)}{|d|^{s+1}}|_{s=0}
-\frac{2\pi}{N}\left(\sum\limits_{a_2=0}^{N-1} \left(\frac{a_2}{N}\right)\right)\sum\limits_{\substack{c=1\\ c\equiv 0(N)}}^{\infty}  \frac{1}{|c|^s}|_{s=0}
=
2L(1,\left(\frac{\cdot}{N}\right)), \]
which further equals $\frac{2\pi h_N}{\sqrt{N}}$. In the first equality we have used $\left(\frac{-1}{N}\right)=-1$, as $N$ is the discriminant of $\QQ[\sqrt{-N}]$. We also compute the non-constant terms:
\[
\frac{-2\pi i}{N}\sum_{a_2=0}^{N-1}\left(\frac{a_2}{N}\right) \sum_{\substack{cn>0, \\ c\equiv 0 (N)}} (\sgn (n)) e^{2\pi i na_2/N} e^{2\pi i nc z/N}
=
\frac{-2\pi i}{N}\sum_{\substack{cn>0 \\ c\equiv 0(N)}}\left(\sum_{a_2=0}^{N-1}\left(\frac{a_2}{N}\right)e^{2\pi i na_2/N}\right)e^{2\pi i nc z/N}.
\]
Computing the Gauss sum $\ds\sum\limits_{a_2=0}^{N-1}\left(\frac{a_2}{N}\right)e^{2\pi i na_2/N}=\left(\frac{n}{N}\right)i\sqrt{N}$,
and changing the notation to $T=nc/N>0$, we get $\frac{4\pi}{\sqrt{N}}\sum\limits_{\substack{T=1}}^{\infty} \sum\limits_{n|T}\sum\limits_{a_2=0}^{N-1}  \left(\frac{n}{N}\right) e^{2\pi i na_2/N} e^{2\pi i T z/N}$.
Finally this gives us the Fourier expansion $\frac{2\pi}{\sqrt{N}}\left(h_N
+
2\sum\limits_{\substack{T=1}}^{\infty} \left(\sum\limits_{n|T, n>0}\left(\frac{n}{N}\right)\right)q^T\right)$ for  $q=e^{2\pi i z}$,
and we recognize that this is $\frac{2\pi}{\sqrt{N}}\Theta_N(z)$.

\bigskip

 Similarly, we also prove:
\begin{Lem}\label{Theta_Eis2}For $N\not\eq 3$, we have $\ds \sum\limits_{c, d \in \ZZ}\frac{\left(\frac{c}{N}\right)}{cz+d}
=
-2 \pi i  \Theta_{N}(z).$ 

\end{Lem}

{\bf Proof:} We write similarly to Lemma \ref{Theta_Eis}, $\ds \sum_{c, d \in \ZZ}\frac{\left(\frac{c}{N}\right)}{cz+d}
=
\sum_{a_1, a_2=0}^{N-1} \left(\frac{a_1}{N}\right) E_1(z; a_1, a_2, N) $. We compute the constant term $\ds \sum_{a_1, a_2=0}^{N-1}\left(\frac{a_1}{N}\right) C(a_1, a_2, N)
=
 - 2\pi  \sum\limits_{c=1}^{\infty} \frac{\left(\frac{c}{N}\right)}{|c|^s}|_{s=0}
=
-2\pi L(0,\left(\frac{\cdot}{N}\right))=-2\pi i h_N. $ For the non-constant part, we get:
\[
-\frac{2\pi i}{N}\sum_{a_1, a_2=0}^{N-1}\left(\frac{a_1}{N}\right) \sum_{\substack{cn>0, \\ c\equiv a_1 (\text{mod } N)}} (\sgn (n)) e^{2\pi i na_2/N} e^{2\pi i nc z/N}
=
-\frac{2\pi i}{N}\sum_{\substack{cn>0}} (\sgn (c)) e^{2\pi i nc z/N}\left(\frac{c}{N}\right)(\sum_{a_2=0}^{N-1}e^{2\pi i na_2/N}), 
\]
As $\sum\limits_{a_2=0}^{N-1}e^{2\pi i na_2/N}\neq 0$ only for $N|c$, we denote $T=nc/N$ and we get $ -4\pi i \sum\limits_{\substack{T=1}}^{\infty} (\sum\limits_{\substack{c|T\\ c>0}}\left(\frac{c}{N}\right)) e^{2\pi i T z},$
which gives us the result.

\bigskip

\subsection{Properties of $\Theta_M$}\label{theta_M}


We will specialize now the previous section to the case $N=3m^*$. Let $M=\QQ[\sqrt{-3m}]$ of discriminant $-3m^*$ and define the Eisenstein series
\begin{equation}\label{eis}
E_{3m^*}(s, z)=\sum_{\substack{c, d \in \ZZ\\ 3m^*|c}}\frac{\left(\frac{d}{3m^*}\right)}{(cz+d)|cz+d|^{s}},
\end{equation}
that will appear in the computation of $L(1, \chi)$ in Section \ref{L}.
For $m\neq 1$, we define the theta function 
\begin{equation}\label{Fourier}
\Theta_M(z)=h_M+2\sum\limits_{n\geq 1} \left(\sum_{d|n} \left(\frac{d}{3m^*}\right)\right)e^{2\pi i n z},
\end{equation}
This is a theta function of weight $1$ for $\Gamma_1(3m^*)$ (see Lemma \ref{modular} below). In general we can define $\Theta_M(z)=h_M+u_M\sum\limits_{n\geq 1} \left(\sum\limits_{d|n} \left(\frac{d}{3m^*}\right)\right)e^{2\pi i n z}$, where $u_M$ is the number of units in $\OO_M$.  We note that for $m=1$ this is the theta function 
\begin{equation}
\Theta_K(z)=1+6\sum\limits_{m, n\in \ZZ}e^{2\pi i (m^2+n^2-mn)z}
\end{equation}
of weight $1$ and level $3$.

The first property we mention is immediate from Lemma \ref{Theta_Eis}:
\begin{cor}\label{SiegelWeil}
$E_{3m^*}(0, z)=\frac{2\pi}{\sqrt{3m^*}}\Theta_M(z).$
\end{cor}

This is a generalization of Theorem 11 in \cite{Ro} that states $E_3(0, z)=\frac{2\pi\sqrt{3}}{9}\Theta_K(z)$.

\bigskip

We present now several properties of $\Theta_M$ that are used in Section \ref{Galois} and Section \ref{X}, in particular the modularity transformation of Lemma \ref{modular} and the inverse transformation of Lemma \ref{transform}.

Using Lemma \ref{Theta_Eis}, we first prove the modular transformation:

\begin{Lem}\label{modular} For $\gamma=\left(\begin{smallmatrix} a & b \\c & d   \end{smallmatrix}\right)\in \Gamma_0(3m^*)$, we have the transformation:
\[
\Theta_M(\gamma z)=\left(\frac{d}{3m^*}\right) (cz+d)\Theta_M(z).
\]
\end{Lem}

{\bf Proof:} Note first that for $\ds \gamma=\left(\begin{smallmatrix} a&b \\ c & d \end{smallmatrix}\right) \in \Gamma_0(3m^*)$, we have
 $E_{3m^*}(s, \gamma z)= \sum\limits_{\substack{p, q \in \ZZ \\ 3m^*|p}} \frac{\left(\frac{q}{3m^*}\right)}{(p\frac{az+b}{cz+d}+q)|p\frac{az+b}{cz+d}+q|^s}$, which can be rewritten as $\left(\frac{d}{3m^*}\right)(cz+d)|cz+d|^s\sum\limits_{\substack{p, q \in \ZZ\\ 3m^*|p}} \frac{\left(\frac{pa+qd}{3m^*}\right)}{((pa+qc)z+pb+qd)|(pa+qc)z+pb+qd|^s}$, or equivalently 
$\left(\frac{d}{3m^*}\right)(cz+d)|cz+d|^s E_{3m^*}(s, z)$.
 Taking the limit $s\ra 0$ as in \cite{Ro}, we get: 
 \[
 E_{3m^*}(0, \gamma z)=\left(\frac{d}{3m^*}\right)(cz+d) E_{3m^*}(0, z)
 \] and together with Lemma \ref{Theta_Eis}, we get the result of the lemma.

\bigskip
We note as a particular case the modular transformation of $\Theta_K$ from \cite{Ro}:
\begin{equation}\label{modular2}
\Theta_K(\gamma z)=\left(\frac{d}{3}\right) (cz+d)\Theta_K(z), \ \  \gamma=\left(\begin{smallmatrix} a&b \\ c & d \end{smallmatrix}\right)  \in \Gamma_0(3),
\end{equation}
which was proved similarly.

Using Lemma \ref{Theta_Eis} and Lemma  \ref{Theta_Eis2}, we show the inverse transformation:
\begin{Lem} We have the transformation
\begin{equation}\label{transform}
\Theta_M(z)=\frac{-1}{3m^*z}\frac{\sqrt{3m^*}}{i}\Theta_M\left(-\frac{1}{3m^*z}\right)
\end{equation}
\end{Lem}

{\bf Proof:} We denote $N=3m^*$ and we have 
\[
\ds \sum_{c, d \in \ZZ}\frac{\left(\frac{d}{N}\right)}{(Ncz+d)|Ncz+d|^s}
= \frac{-1}{Nz|Nz|^{s}}\sum_{c, d \in \ZZ}\frac{\left(\frac{d}{N}\right)}{(c+d\frac{-1}{Nz})|c+d\frac{-1}{Nz}|^s}.
\] Taking the limit to $s=0$, we get from Lemma \ref{Theta_Eis} and Lemma \ref{Theta_Eis2}:
\[
\frac{2\pi }{\sqrt{3m^*}}\Theta_M(z)=-2\pi i  \frac{-1}{3m^*z} \Theta\left(\frac{-1}{3m^*z}\right),
\]
which gives us the result.

\bigskip
This is a generalization of the transformation for $m=1$ mentioned in \cite{Ro}:
\begin{equation}\label{transform2}
\Theta_K(z)=\frac{-1}{3z}\frac{\sqrt{3}}{i}\Theta_K\left(-\frac{1}{3z}\right).
\end{equation}
\bigskip

Finally, we present a result relating sums of theta functions that will be used in Section \ref{X}:

\begin{Lem}\label{sum_M} $\ds\sum_{0\leq j \leq d-1}\Theta_M\left(z+\frac{j}{d}\right)=d\Theta_M(dz)$ for $d|3m^*$.
\end{Lem}

{\bf Proof:} Using the Fourier expansion \eqref{Fourier} $\Theta_M(z)=\sum\limits_{N\geq 0} c_N e^{2\pi i N z}$, with $c_0=1$, $c_N=2\sum\limits_{d|N} \left(\frac{d}{3m^*}\right)$ for $N\geq 1$, we rewrite:
\[
\sum_{0\leq j\leq d-1}
\Theta_M\left(z+\frac{j}{d}\right)
=
\sum_{N\geq 0} c_N \sum_{0\leq j\leq d-1}
e^{2\pi i N (z+j/d)}
=
\sum_{N\geq 0} c_N e^{2\pi i N z}\sum_{j=0}^{d-1}e^{2\pi i N j/d},
\]
which equals $ d\sum\limits_{\substack{N\geq 0\\ d|N}}c_{N} e^{2\pi i Nz}.$ As $c_N=c_{N/d}$, we get on the LHS the Fourier expansion $d\sum\limits_{\substack{N'\geq 0}}c_{N'} e^{2\pi i N'dz}=d\Theta_{M}(dz)$.

\section{Computing $L(1, \chi)$}\label{L}

In this section we use the method of Rosu \cite{Ro} to compute the value of $L(1, \chi)$ using Tate's thesis as a linear combination of values of theta functions and Hecke characters.

\subsection{Schwartz-Bruhat functions}\label{Schwartz}
We choose the Schwartz-Bruhat function $\Phi_f\in \SSS(\AAA_{K, f})$ such that Tate's zeta function $Z(s, \Phi, \chi)$ defined below to be nonzero. More precisely, $\Phi_f=\prod\limits_{v\nmid \infty} \Phi_v$, where:

\begin{itemize}
\item $\Phi_v=\Char_{\OO_{K_v}}$ for $v\nmid 3D$,
\item $\Phi_p=\sum\limits_{(a, p)=1}\left(\frac{a}{p}\right)\Char_{(a+p\OO_{K_p})}$ for $p|m$,

\item $\Phi_p=\Char_{(\ZZ+D\OO_{K_p})}$ for $p|D$

\item $\Phi_v=\Char_{(1+3\OO_{K_v})}$ for $v=\sqrt{-3}$.

\item $\Phi_2=\begin{cases} \Char_{\ZZ_2[\omega]} & \text{ for } \alpha\equiv 1(4) \\
\Char_{\pm 1+4\ZZ_2[\omega]} & \text{ for } \alpha\equiv -1(4) \\
 \Char_{1+4\ZZ_2[\omega]}-\Char_{3+4\ZZ_2[\omega]}  & \text{ for }  m\equiv 3 (4).\\
\end{cases}$
\end{itemize}

We note that we have chosen the local Schwartz-Bruhat functions such that 
\begin{equation}\label{phi_chi}
\begin{cases} \Phi_v\chi_v=\Char_{\ZZ+3m^*D'\OO_{K_v}} & v|m^*D' 
\\ 
\Phi_{v}\chi_{v}=\Char_{1+3\OO_{K_v}} & v=\sqrt{-3}, \end{cases}
\end{equation}
 where $\chi_v$ the local component of the Hecke character $\chi:\AAA_K^{\times}/K^{\times}\ra \CC$ defined in \eqref{chi}, corresponding to the elliptic curve $E_{D, \alpha}$ from CM theory. 
 
 We also recall $D'=\begin{cases}
D & \text{for } \alpha \equiv 1, \pm \sqrt{-3} (4) \\ 4D &  \text{for } \alpha \equiv -1(4)\end{cases}$.

\subsection{Tate's zeta function} For a Hecke character $\chi:\AAA_{K}^{\times}/K^{\times} \ra \CC^{\times}$ and a Schwartz-Bruhat function $\Phi\in \SSS(\AAA_{K})$, Tate's zeta function is defined locally as $\ds Z_v(s, \chi_v, \Phi_v)=\int\limits_{K_v^{\times}} \chi_{v}(\alpha_v) |\alpha_v|^s_v \Phi_v(\alpha_v) d^{\times} \alpha_v$, where we choose the self-dual Haar measure as in Tate's thesis (see \cite{Bu}). Globally $\ds Z(s, \chi, \Phi)=\prod\limits_v Z_v(s, \chi_v, \Phi_v)$ has meromorphic continuation to all $s\in \CC$ and in our case it is entire.

We will compute $Z_f(s, \chi_f, \Phi_f)=\prod\limits_{v\nmid \infty} Z_v(s, \chi_v, \Phi_v)$ for $\chi$ the Hecke character \eqref{chi} corresponding to $E_{D, \alpha}$ via CM theory, and the Schwartz-Bruhat function $\Phi_f$ chosen above in Section \ref{Schwartz}. From Tate's thesis (see \cite{Bu}, Proposition 3.1.4), we have the equality of local factors $L_v(s, \chi)=Z_v(s, \chi)$ at all the unramified places $v\nmid 3m^*D'$. We define for $v|\overline{\alpha}$ the local $L$-function:
\[
L_{\overline{\alpha}}(s, \chi)=\prod_{v|\overline{\alpha}} L_{v}(s, \chi),
\] 
where $L_{v}(s, \chi)=(1-\eps(\pp_v)\chi_D(\pp_v)\varpi_vq_v^{-s})^{-1}$ for 
$\varpi_v$ the unique generator of $\pp_v$ such that $\varpi_v\equiv 1 (3)$, and $q_v=\Nm\alpha_v$.

We fix the Schwartz-Bruhat function $\Phi_f$ from section \ref{Schwartz} and it will be an immediate consequence of Tate's thesis:
\begin{Lem} For all $s\in \CC$, we have:
\begin{equation}\label{L_to_Z} 
\frac{L(s, \chi)}{L_{\overline{\alpha}}(1, \chi)}=Z_f(s, \chi_f, \Phi_f)\frac{1}{2}\prod_{p|3D'm^*}\vol(\left(\ZZ+3D'm^*\OO_{K_p}\right)^{\times}).
\end{equation}

\end{Lem}

{\bf Proof:} As $L_v(s, \chi)=Z_v(s, \chi)$ at all the unramified places $v\nmid 3m^*D'$, we have the equality $\ds L(s, \chi)=Z_f(s, \chi, \Phi)\prod\limits_{v|3D'm^*}\frac{L_v(s, \chi)}{Z_v(s, \chi, \Phi_v)}.$ At the ramified places of the $L$-function, we have $\Phi_v(x)\chi(x)=1$, when $\Phi_{v}$ is nonzero. Thus, from \eqref{phi_chi}, we have $\prod\limits_{v|3m^*D'} Z_v(s, \chi, \Phi_p)=\frac{1}{2}\prod\limits_{p|3D'm^*}\vol(\left(\ZZ+3D'm^*\OO_{K_p}\right)^{\times}).$ The $1/2$ factor occurs from considering $\vol(1+3\OO_{K_{v}})$ at $v=\sqrt{-3}$. Finally, by definition, the terms $L_v(s, \chi)=1$ for $v|3D\alpha$ for $\alpha\equiv 1(4)$, respectively for $v|6D\alpha$ for $\alpha\equiv -1, \pm \sqrt{-3}(4)$. Thus $\prod\limits_{v|3m^*D'}L_v(s, \chi)=L_{\overline{\alpha}}(s, \chi)$.

\bigskip

Next we use \eqref{L_to_Z} to get the value of $L(s, \chi)$ by computing the value of $Z_f(s, \chi_{f}, \Phi_f)$ as a linear combination of Hecke characters:

\begin{Lem}\label{Z_linear} For all $s\in \CC$, we have:
\[
\frac{L(s, \chi)}{L_{\overline{\alpha}}(s, \chi)}=\sum_{\beta_f \in U(3D'm^*)\setminus\AAA_{K, f}^{\times}/K^{\times}} I(s, \beta_f, \Phi_f)\chi_f(\beta_f),
\]
\noindent where $\ds I(s, \beta_f, \Phi_f)=\sum\limits_{k\in K^{\times}} \frac{k}{|k|_{\CC}^{2s}} \Phi_f(k\beta_f)$ and $U(3D'm^*)=(1+3\OO_{K_{\sqrt{-3}}})\prod\limits_{p|D'm^*}{(\ZZ+D'm^*\OO_{K_p})^{\times}} \prod\limits_{v\nmid 3D}{\OO_{K_v}^{\times}}$.

\end{Lem}

{\bf Proof:} We take the quotient by $K^{\times}$ in the integral defining $Z_f(s, \chi_f, \Phi_f)$. Noting that $\chi_{f}(k\beta_f)=\chi_{\infty}^{-1}(k)\chi_{f}(\beta_f)=k\chi_f(\beta_f)$ and $|k\beta_f|^s_f=|k|_{\CC}^{-2s}|\beta_f|_f^s$, where $| \cdot |_{\CC}$ is the usual absolute value over $\CC$, we get:
\[
Z_f(s, \chi_f, \Phi_f)
=
\int\limits_{\AAA_{K, f}^{\times}/K^{\times}}\left(\sum_{k\in K^{\times}} \frac{k}{|k|_{\CC}^{2s}} \Phi_f(k\beta_f)\right) \chi_{f}(\beta_f)|\beta_f|_f^s ~d^{\times} \beta_f.
\]

As $\varphi$ and $|\cdot|$ are unramified on $U(3m^*D')$, and $\Phi_f(k\beta_f)\chi_{f}(\beta_f)=1$ on $U(3m^*D')$, we get a finite sum $Z_f(s, \chi_D\varphi, \Phi_f)
=
\vol(U(3D'm^*))\sum \limits_{U(3D'm^*)\setminus \AAA_{K, f}^{\times}/K^{\times}}I(s, \beta_f, \Phi_f) \chi_{f}(\beta_f) |\beta_f|_f^s.$ As  $\vol(U(3D'm^*))=\frac{1}{2}\prod\limits_{p|3D'm^*}\vol(\left(\ZZ+3D'm^*\OO_{K_p}\right)^{\times})$ and using \eqref{L_to_Z}, we get the result of the lemma.

\bigskip

We recall the Eisenstein series from Section \ref{theta_M}:
\begin{equation}\label{eis}
E_{3m^*}(s, z)=\sum_{\substack{c, d \in \ZZ\\ 3m^*|c}}\frac{\left(\frac{d}{3m^*}\right)}{(cz+d)|cz+d|^{s}}.
\end{equation}

Then we have:

\begin{Lem}\label{lin_comb1}  $\ds 
\frac{L(s, \chi)}{L_{\overline{\alpha}}(s, \chi)}
=
\sum\limits_{[\A]\in \Cl(\OO_{3m^*D'})} E_{3m^*}(2s-2, D'\tau_{\A})
\chi(\A)a^{1-2s}.$

\end{Lem}

{\bf Proof:} For representatives $\beta_f \in \prod\limits_{v} \OO_{K_v}^{\times}$, $\beta_3\equiv 1(3)$, we take $\beta_0\in \OO_K$ such that $\beta_0\equiv \beta \mod 3Dm^*$. Thus, for $k\in K$, we have $\Phi_f(k\beta_f)=\Phi_f(k\beta_0)$ and this is nonzero only for $k\in \OO_K$.  Let $\A=(\beta_0)$ be the corresponding ideal.   Moreover, $\varphi_{f}(\beta_f)=\varphi_{f}((\beta_0)_{p|3D})=1$, $\chi_{D, f}(\beta_f)=\chi_{D, f}((\beta_0)_{p|3D})=\chi_{D, v}^{-1}((\beta_0)_{v|\beta_0})=\chi_D(\overline{\A})$ and similarly $\eps_{f}(\beta_f)=\eps^{-1}(\A)=\eps(\A)$. We note $U(3m^*D')\setminus \AAA_{K, f}^{\times}/K^{\times}\simeq \Cl(\OO_{3m^*D'})$, the ring class field corresponding to the order $\OO_{3m^*D'}=\ZZ+3m^*D'\OO_K$.  Thus we have:
\[
\frac{L(s, \chi)}{L_{\overline{\alpha}}(s, \alpha)}=
\sum\limits_{[\A] \in \Cl(\OO_{3m^*D'})} \left(\sum_{k\in K^{\times}} \frac{k}{|k|_{\CC}^{2s}}  \Phi_f(k\beta_0)\right) \chi_{D}(\overline{\A})\eps(\A), 
\]
where the ideals $\A$ are representatives of the classes of $\Cl(\OO_{3m^*D'})$. 

As $k\in \OO_K$, we have $\Phi_{3Dm^*}(k\beta_0)\neq 0 $ for $k\beta_0 \in \A$ and $k\beta_0 \in (\ZZ+3Dm^*\OO_{K})^{\times}$ with $k\beta_0 \equiv 1(3)$. Write the ideal $\A=\left[a, \frac{-b+\sqrt{-3}}{2}\right]$ with $a=\Nm\A$, $b^2\equiv -3(4a)$. Thus:
\[
\frac{|\beta_0|_{\CC}^{2s}}{\beta_0}\sum_{k\in \OO_{K}} \frac{\beta_0 k}{|\beta_0 k|_{\CC}^{2s}} \Phi_{3m^*D'}(\beta_0 k)
=
\frac{a^{s}}{\varphi(\A)}\sum_{\substack{c, d \in \ZZ \\ d\equiv 1(3) }} \frac{\left(\frac{da}{m^*}\right)(da+3Dm^* c\frac{-b+\sqrt{-3}}{2})}{
|da+3m^*D' c\frac{-b+\sqrt{-3}}{2}|^{2s}}.
\]
We further rewrite this as $\frac{1}{2}\frac{\overline{\varphi(\A)}}{a^s}\left(\frac{a}{3m^*}\right)E_{3m^*}(2s-2, D'\tau_{\overline{\A}})$.
Note further that $\left(\frac{a}{3}\right)=1$, $\left(\frac{a}{m^*}\right)=\left(\frac{m}{a}\right)=\left[\frac{m}{\A}\right]=\left[\frac{\alpha}{\A}\right]\left[\frac{\overline{\alpha}}{\A}\right]=\eps(\A)\eps(\overline{\A})$, and we get the sum: 
\[
\frac{L(s, \chi)}{L_{\overline{\alpha}}(s, \chi)
}
=
\frac{1}{2}
\sum\limits_{[\A]\in \Cl(\OO_{3m^*D'})} E_{3m^*}(2s-2, D'\tau_{\overline{\A}})
\chi(\overline{\A}) a^{1-2s}.
\]
Changing $\A$ to $\overline{\A}$ in  $\Cl(\OO_{3m^*D'})$ we get the result of the lemma.

\bigskip

We recall the theta function $\Theta_K(z)=1+6\sum\limits_{m, n\in \ZZ}e^{2\pi i (m^2+n^2-mn)z}$ of weight $1$ and level $3$. From Lemma 3.8 of \cite{Ro}, for $\A=\left[a, \frac{-b+\sqrt{-3}}{2}\right]$, we have:
\begin{equation}\label{lemma_phi}
\Theta_K\left(\tau_{\A}\right)
=
\varphi(\overline{\A})\Theta_K(\omega)
\end{equation}

Plugging in $s=1$ in Lemma \ref{lin_comb1}, together with \eqref{lemma_phi}, we get:

\begin{Prop}\label{lin_comb}
$L(1, \chi)
=
\frac{L_{\overline{\alpha}}(1, \chi)\Theta_K(\omega)}{2D^{1/3}\alpha^{1/2}}
\sum\limits_{[\A]\in \Cl(\OO_{3m^*D'})} \frac{E_{3m^*}(0, D'\tau_{\A})}{\Theta_{K}(\tau_{\A})}
\chi_{D}(\A)\eps(\A) D^{1/3} \alpha^{1/2}.$
\end{Prop}

We recall from Corollary \ref{SiegelWeil} that we have:
\begin{equation}
E_{3m^*}(0, z)=\frac{2\pi}{\sqrt{3m^*}}\Theta_M(z).
\end{equation}

Thus we get in Proposition \ref{lin_comb}:
\[
L(1, \chi)
=
L_{\overline{\alpha}}(1, \chi) \frac{\pi }{\sqrt{3m^*} }\frac{\Theta_K(\omega)}{D^{1/3}\alpha^{1/2}}
\sum\limits_{[\A]\in \Cl(\OO_{3m^*D'})} \frac{\Theta_M(D'\tau_{\A})}{\Theta_{K}(\tau_{\A})}
\chi_{D}(\A)\left[\frac{\alpha}{\A}\right] D^{1/3} \alpha^{1/2}.
\]

Recall from \cite{Ro} that $\Theta_K\left(\omega\right)=\Gamma\left(\frac{1}{3}\right)^3/(2\pi^2)$. Denoting:

\begin{equation}\label{T_D}
X_{D\alpha}=\sum\limits_{[\A]\in \Cl(\OO_{3m^*D'})} \frac{\Theta_M(D'\tau_{\A})}{\Theta_{K}(\tau_{\A})}
\chi_{D}(\A)\left[\frac{\alpha}{\A}\right] D^{1/3} \alpha^{1/2},
\end{equation}
we get:

\begin{cor}\label{form_1} $\ds L(E_{D\alpha/K}, 1)=\Omega_{D\alpha}\overline{\Omega}_{D\alpha} \frac{4|L_{\overline{\alpha}}(1, \chi)|^2}{3m^*} |X_{D\alpha}|^2$.
\end{cor}

\section{Galois conjugates}\label{Galois}

Our final step in showing Theorem \ref{thm1} is to prove that the individual terms in the sum $X_{D\alpha}$ in \eqref{T_D} are Galois conjugates to each other, using Shimura's reciprocity law. For $f$ a modular function of level $N$ and $\tau \in \HH$ a CM point, $f(\tau)$ is an algebraic integer and Shimura's reciprocity law gives an explicit way to compute its Galois conjugates. More precisely, for $\tau\in \HH\cap K$ and $\sigma\in \Gal(K^{ab}/K)$ Shimura's reciprocity law states:
\[
f(\tau)^{\sigma} =f^{g_{\sigma}}(\tau),
\]
where $g_{\sigma}\in \GL_2(\AAA)$ acts on the modular function $f$. For the explicit action see \cite{La}, \cite{Ro} or \cite{G}.

We will recall an explicit version of Shimura's reciprocity law. For a general treatment see \cite{St}, \cite{GS}, \cite{G}, as well as \cite{Ro} for more details for the current approach. Let $f$ be a modular function of level dividing $3N$ and $\A$ an primitive ideal prime to $3N$. We write $\A=[a, \frac{-b+\sqrt{-3}}{2}]$ and let $k_{\A}=ta+s\frac{-b+\sqrt{-3}}{2}$ be a generator of $\A$ and $c=\frac{b^2+3}{4a}$. Then for $\sigma_{\A}$ the Galois action corresponding to the ideal $\A$ via the Artin map, we have for $\tau=\frac{-b+\sqrt{-3}}{2}$:
\[
(f(\tau))^{\sigma_{\A}^{-1}}=f^{\left(\begin{smallmatrix} ta-sb & -sca \\s & ta   \end{smallmatrix}\right)_{p|3N}}(\tau),
\]
Moreover, if $f$ has rational coefficients at $\infty$, we have:

\begin{equation}\label{Shimura}
(f(\tau))^{\sigma_{\A}^{-1}}=f^{\left(\begin{smallmatrix} ta-sb & -sc \\s & t   \end{smallmatrix}\right)_{p|3Dm^*}}(\tau)=f(\left(\begin{smallmatrix} ta-sb & -sc \\s & t   \end{smallmatrix}\right)\tau).
\end{equation}

We also recall a different explicit version of Shimura's reciprocity law, that is Lemma 5.4 from \cite{Ro}. For a modular function $f$ of level dividing $3N$ and a primitive ideal $\A=[a, \frac{-b+\sqrt{-3}}{2}]$ prime to $3N$, we have the Galois action:
\begin{equation}\label{Shimura_A}
f(\tau)^{\sigma_{\A}^{-1}}=f(\tau_{\A}),
\end{equation}
where $\tau_{\A}=\frac{-b+\sqrt{-3}}{2a}$. 

We will use \eqref{Shimura} and \eqref{Shimura_A} at various points in the paper.  At the moment we are interested in the Galois conjugates of the modular function $F(z)=\frac{\Theta_M(D'z)}{\Theta_K(z)}$ at $z=\tau$. Our goal is to show that $F(\tau)D^{1/3}\alpha^{1/2}\in H_{3D'm^*}$ and its Galois conjugates are the terms appearing in the sum defining $X_{D\alpha}$.

Using Shimura's reciprocity law \eqref{Shimura} we compute the Galois conjugates of $F(\tau)$ under $\Gal(K^{ab}/H_{3D'm^*})$ below:

\begin{Lem} For  $\A=[a, \frac{-b+\sqrt{-3}}{2}]$ a primitive ideal with generator $k_{\A}=ta+s\frac{-b+\sqrt{-3}}{2}$, $3m^*D'|s$, $t\equiv 1(3)$, we have:
\[
(F(\tau))^{\sigma_{\A}^{-1}}=\left(\frac{t}{m^*}\right)F(\tau).
\]

\end{Lem}

{\bf Proof:} Using the remarks above, for  $c=\frac{b^2+3}{4a}$, and $\tau=\frac{-b+\sqrt{-3}}{2}$ we need to compute:
\[
F(\left(\begin{smallmatrix} ta-sb & -sc \\s & t   \end{smallmatrix}\right)\tau)=
\frac{\Theta_M(\left(\begin{smallmatrix} ta-sb & -scD' \\s/D' & t   \end{smallmatrix}\right)D'\tau)}{\Theta_K(\left(\begin{smallmatrix} ta-sb & -sc \\ s & t   \end{smallmatrix}\right)\tau)}
\]
Using the modular transformations proved in Lemma \ref{modular} for $3m^*|(s/D')$ and \eqref{modular2} from Section \ref{theta_M}, we get 
\[
F(\left(\begin{smallmatrix} ta-sb & -sc \\s & t   \end{smallmatrix}\right)\tau)=
\left(\frac{t}{m^*}\right)\frac{\Theta_M(D'\tau)}{\Theta_K(\tau)},
\]
which proves the lemma.

\bigskip
It follows that

\begin{Lem}\label{H_3mD} $F(\omega)\alpha^{1/2}D^{1/3}\in H_{3m^*D'}$.
\end{Lem}

{\bf Proof:} We want to show that $F(\omega)\alpha^{1/2}D^{1/3}$ is invariant under the action of $\Gal(H_{3m^*D'}/K)$. Let $\A=[a, \frac{-b+\sqrt{-3}}{2}]$ a primitive ideal with generator $k_{\A}=ta+s\frac{-b+\sqrt{-3}}{2}$, $3m^*D'|s$, $t\equiv 1(3)$. Then $(\alpha^{1/2})^{\sigma_{\A}^{-1}}=\left[\frac{\alpha}{k_{A}}\right]\alpha^{1/2}$ and by quadratic reciprocity (see Section \ref{quad_char},  \ref{rec_law0}), we have $\left[\frac{\alpha}{k_{A}}\right]=\left[\frac{k_{A}}{\alpha}\right]$, for $\alpha \equiv \pm 1 (4)$ and $\left[\frac{\alpha}{k_{A}}\right]=\left[\frac{k_{A}}{\alpha}\right]\left(\frac{ta}{4}\right)$ for $\alpha \equiv \pm \sqrt{-3} (4)$. Thus we get 
$\left[\frac{\alpha}{k_{A}}\right]=\left[\frac{ta}{\alpha}\right]=\left(\frac{ta}{m}\right)$ for $m\equiv 1 (4)$ and $\left[\frac{\alpha}{k_{A}}\right]=\left[\frac{ta}{4\alpha}\right]=\left(\frac{ta}{4m}\right)$ for $m\equiv 3(4)$. Combining this with the lemma above, we get:
\[
(F(\omega)\alpha^{1/2})^{\sigma_{\A}^{-1}}=\left(\frac{a}{m^*}\right)F(\omega)\alpha^{1/2}.
\]

Note that as $3m^*|s$, we have $a=\Nm k_{\A} \equiv t^2a^2 \mod 3m^*$, thus $1/t^2\equiv a \mod 3m^*$, implying $\left(\frac{a}{m^*}\right)=1$. As $D^{1/3}\in H_{3D'm^*}$, this finishes our proof.

\bigskip

Moreover, using \eqref{Shimura_A} for the modular function $F(z)=\frac{\Theta_M(Dz)}{\Theta_K(z)}$, we have for $\tau_{\A}=\frac{-b+\sqrt{-3}}{2a}$:
\[
F(\tau)^{\sigma_{\A}^{-1}}=F(\tau_{\A}).
\]
Furthermore, by definition we have $(D^{1/3})^{\sigma_{\A}^{-1}}=\chi_D(\A)D^{1/3}$ and $(\alpha^{1/2})^{\sigma_{\A}^{-1}}=\left[\frac{\alpha}{\A}\right]\alpha^{1/2}$. Then we can rewrite \eqref{T_D}:
\begin{equation}\label{T_lin_comb}
X_{D\alpha}
=
\sum\limits_{[\A]\in \Cl(\OO_{3m^*D'})} (F(\tau)D^{1/3} \alpha^{1/2})^{\sigma_{\A}^{-1}}.
\end{equation}

\bigskip

Then combining with Lemma \ref{H_3mD}, we get:

\begin{Prop} $X_{D\alpha}=\Tr_{H_{3D'm^*}/K} (F(\omega)D^{1/3} \alpha^{1/2})$ and $X_{D\alpha}\in \OO_K$.
\end{Prop}

To see that $X_{D\alpha}\in \OO_K$, it is enough to show that $F(\omega)$ is an algebraic integer. This follows from the fact that $\omega$ is a CM point and $F(\gamma z)$ is a modular function that has integer coefficients in its Fourier expansion at $\infty$ for all $\gamma \in \SL_2(\ZZ)$.

\bigskip
Together with Corollary \ref{form_1}, this implies:

\begin{Thm}\label{value} $\ds 
L(E_{D\alpha/K}, 1)
=
\Omega_{D\alpha}\overline{\Omega}_{D\alpha} c \left|\Tr_{H_{3D'm^*}/K} \left(\frac{\Theta_M(D'\omega)}{\Theta_K(\omega)}D^{1/3} \alpha^{1/2}\right)\right|^2,$ where $c=\frac{4|L_{\overline{\alpha}}(1, \chi)|^2}{3m^*}$.

\end{Thm}

\section{Value as a square}\label{X}
To simplify the calculations, we only consider the case of $\alpha$ prime. The goal of this section is to show Theorem \ref{thm2} from the Introduction, that we restate below as Theorem \ref{thm_square}. Starting with Theorem \ref{value} that states
\[
L(E_{D\alpha/K}, 1)
=
\frac{4^e}{3}\Omega_{D\alpha}\overline{\Omega}_{D\alpha} c \left|\frac{L_{\overline{\alpha}}(1, \chi)}{\alpha}X_{D\alpha}\right|^2, \ \ e=\begin{cases} 1 & m\equiv 1(4) \\ 0 & m\equiv 3(4) \end{cases},
\] 
we will show that $L_{\overline{\alpha}}(1, \chi)X_{D\alpha}/\alpha$ equals, up to a sextic root of unity, the following trace:
\begin{equation}
 Z=\omega^l\Tr_{H_{3D^*}/K} 
 \left(\frac{\Theta_M(D^*\tau_0/m^*)}{\Theta_K(\tau_0/m)}D^{1/3} \alpha^{1/2}\right),
\end{equation}
where $\tau_0=\frac{B+\sqrt{-3}}{2}$ with $B \equiv \sqrt{-3}(\alpha)$, $B$ odd and $D^*=\begin{cases} D & \alpha\equiv 1(4) \\ 4D & \alpha\equiv -1, \pm\sqrt{3}(4)\end{cases}$. 



Using Theorem \ref{value}, we will show:

\begin{Thm}\label{thm_square} Let $\alpha\in \OO_K$ prime and $D$ any integer such that $(D, 6\alpha)=1$, $(\alpha, 3)=1$. Then:
\[ 
\frac{L(E_{D\alpha/K}, 1)}{\Omega_{D\alpha}\overline{\Omega}_{D\alpha}}=c_0(2^eZ_{D\alpha}/\sqrt{-3})^2, 
\]
where $c_0=(-1)^{v_3(c_{E_{D, \alpha}})}$ and $e=1$ for $m\equiv 1(4)$, $e=0$ for $m\equiv 3(4)$ as above.

Moreover, $Z\in \OO_K$ and $\overline{Z/\sqrt{-3}}=c_0(Z/\sqrt{-3})$. In particular, $Z/\sqrt{-3} \in \ZZ$ for $c_0=1$ and $Z\in \ZZ$ for $c_0=-1$, respectively.

\end{Thm}


Denote $S_{D\alpha}=\frac{L(E_{D\alpha/K}, 1)}{\Omega_{D\alpha}\overline{\Omega}_{D\alpha}c_{E_{D, \alpha}}}(\#E_{D, \alpha}(K)_{tor})^2$. As $\#E_{D, \alpha}(K)_{tor}=1$ for $\Nm\alpha>1$, we will immediately get that $S_{D\alpha}$ is a square when 	
\begin{itemize}
		\item $c_{E_{D, \alpha}}=1, 4$ for $m\equiv 1(4)$
	\item $c_{E_{D, \alpha}}=1$ for $m\equiv 3(4)$
	\end{itemize}
	

These are the situations described below:

\begin{cor}\label{square_sha} For $m\equiv 1(4)$, $\Sha_{an, E_{D, \alpha}}$ is an integer square in the cases:

\begin{itemize}
	
	\item $D\equiv \pm 1(9)$, $\alpha\equiv -1(3)$ and $\left[\frac{\alpha}{\pp}\right]=-1$ for all $\pp|D$
	
	\item $D\equiv \pm 4(9)$, $\alpha\equiv 1(3)$ and $\left[\frac{\alpha}{\pp}\right]=-1$ for all $\pp|D$
	
	\item $D\equiv \pm 2(9)$, $\chi_{2D^2}(\alpha)\neq1$
\end{itemize}
	
	 For $m\equiv 3(4)$, $\Sha_{an, E_{D, \alpha}}$ is an integer square in the cases:
\begin{itemize}
	\item $D\equiv \pm 1(9)$ and $\left[\frac{\alpha}{\pp}\right]=-1$ for all $\pp|D$, $\chi_{2D^2}(\alpha)\neq 1$
	
	\item $D\equiv \pm 4(9)$, $\alpha\equiv 1(3)$ and $\left[\frac{\alpha}{\pp}\right]=-1$ for all $\pp|D$, $\chi_{2D^2}(\alpha)\neq 1$

\end{itemize}

\end{cor}

For $D=1$ and $\Nm\alpha>1$, we have $(E_{\alpha}(K))_{tor}=\{\infty\}$ and $c_{E_{\alpha}}\in \{1, 3\}$, thus we get:

\begin{cor} For $\alpha\equiv -1(3)$ prime, $S_{\alpha}$ is an integer square. 
\end{cor}


In order to show Theorem \ref{thm_square}, we want to change $L_{\overline{\alpha}}(1, \chi) X/\alpha$ by a cubic root of unity and show that this is an integer or $\sqrt{-3}$ times an integer. The proof proceeds as follows. We introduce several traces $T, U, W_i, Y_i$ and $Y$ (see the beginning of Section \ref{XT} and Section \ref{XYZ}).

In Section \ref{XT} we compute two formulas that relate $X$ to $T$ and $U$, showing in Corollary \ref{real} that:
\[
L_{\overline{\alpha}}(1, \chi)X
=
U,
\]
with $U/\alpha$ by definition being equal to $Z$ up to a cubic root of unity.
The main goal of section \ref{XYZ} is to relate $U/\alpha$ to its complex conjugate and we start by noting that $U$ equals $Y/D'$ up to a sextic root of unity. Thus we are interested in the complex conjugate of $Y/\alpha$. We first compute several relations between $Y_j$ and $Y$, and their complex conjugates,
which leads to showing in Proposition \ref{conjY} and Proposition \ref{conjU} that $Y/\alpha$, and respectively $U/\alpha$, are real or purely imaginary up to a third root of unity.

This culminates with the proof of Corollary \ref{real} in which we write up $\frac{L(E_{D, \alpha}, 1)}{\Omega_{D\alpha}\overline{\Omega}_{D\alpha}}$ as $4^e|Z|^2/3$, where $Z$ is real or purely imaginary, with $Z$ equal to $U/\alpha$ up to a cubic root of unity. We finally show in Corollary \ref{integral} that $Z$ is an integer or $\sqrt{-3}$ times an integer, which finishes the proof of Theorem \ref{thm_square}.


\bigskip

We note that we have to treat the cases $m\equiv 1(4)$ and $m\equiv 3(4)$ separately, due to the formula for the transformation \eqref{transform} of $\Theta_M$ from $z$ to $-1/3m^*z$. Thus we will define different traces for $m\equiv 1(4)$ and $m\equiv 3(4)$ with the change:
\[
\Tr_{H_{3D'2^{e_0}m}}\frac{\Theta_M(D'^{f}\frac{\tau/2^{e_0}}{m^{e'}3^{e''}})}{\Theta_K(\frac{\tau}{3^{e'''}})}, \text{where } e_0=\begin{cases}
0 & m\equiv 1(4) \\
1 & m\equiv 3(4)
\\\end{cases}.
\]
Here $f=\pm 1$, $e', e'', e'''\in \{0, 1\}$.

Throughout the proofs we will use several results proved in Section \ref{comp_shimura} that are applications of the Shimura reciprocity law. In particular, we show that the traces we take are well defined and take values in $\OO_K$. Moreover, the computations of the characters $\chi_D$ and $\eps$ at various ideals, as well as the explicit computation of various Galois conjugates of modular functions are treated in the Appendix.

\subsection{Computing $X$}\label{XT}

For $m\equiv 3(4)$, we define $X^{(i)}=\Tr_{H_{6Dm}/K} \frac{\Theta_M\left(D\tau_i/2\right)}{\Theta_K (\omega)} \alpha^{1/2}D^{1/3}$ with $i\in \{1, 3\}$, $b_1\equiv 1(8)$, $b_3\equiv -1(8)$. Then from Lemma \ref{X_1} in Section \ref{comp_shimura}, we have $X=X^{(i)}$ for $m\equiv 3(4)$ under the condition 
\begin{equation}\label{cond}
[\alpha]\left(\frac{-b_i}{4}\right)=1,
\end{equation}
where the symbol $[\alpha]$ was defined in Section \ref{quad_char}.

We will fix $b_{i}(8)$ to satisfy \eqref{cond} and we will use the notation $X^*$ for $X^{(i)}$ when $m\equiv 3(4)$ and for $X$ when $m\equiv 1(4)$. Clearly $X=X^*$ and our goal will be to compute $L_{\overline{\alpha}}(1, \chi) X^*/\alpha$, for which we get:
\[
\frac{L(E_{D, \alpha}, 1)}{\Omega_{D\alpha}\overline{\Omega}_{D\alpha}}
=
\frac{4}{3m^*}|L_{\overline{\alpha}}(1, \chi){X^*}^2|.
\]

In this section, we define the traces $U$ and $T$(see table below) and compute two formulas that relate $X^*$ to $T$ and to $U$ in Lemmas \ref{rel_1} and \ref{rel_2}. They give us Proposition \ref{rel_3}, the main result of the section, in which we show that $L_{\overline{\alpha}}(1, \chi)X=U$.

\bigskip
Below are the important traces we take:

\begin{center}
\begin{tabular}{|c |l |l |}
\hline
& $\alpha \equiv \pm 1(4)$ &$\alpha \equiv \pm \sqrt{-3}(4)$ \\ 
\hline
$X^*$ & $\Tr_{H_{3D'm}/K} \frac{\Theta_M\left(D'\omega\right)}{\Theta_K (\omega)} \alpha^{1/2}D^{1/3}$ &$\Tr_{H_{6Dm}/K} \frac{\Theta_M\left(D\tau_1/2\right)}{\Theta_K (\omega)} \alpha^{1/2}D^{1/3}$  \\
\hline
$\ds U$ &
$ 
\ds \Tr_{H_{3D'}/K}
\frac{\Theta_M\left(D'\tau_{0}^*/m\right)}{\Theta_K\left(\omega\right)} \alpha^{1/2}D^{1/3}$ &
$ 
\ds \Tr_{H_{6D}/K}
\frac{\Theta_M\left(D\frac{\tau_{0}^*}{2m}\right)}{\Theta_K\left(\omega\right)} \alpha^{1/2}D^{1/3}$ 
\\
\hline
$\ds T $
&
$
 \Tr_{H_{3D'm}/K}
\frac{\Theta_M\left(\frac{-\overline{\tau}_1}{3D'}\right)}{\Theta_M\left(-\overline{\tau}_1/3\right)} \overline{\alpha}^{1/2}D^{1/3}$ & 
$ \Tr_{H_{6Dm}/K}
\frac{\Theta_M\left(\frac{-\overline{\tau}_1}{6D}\right)}{\Theta_M\left(-\overline{\tau}_1/3\right)} \overline{\alpha}^{1/2}D^{1/3}$  
\\ \hline
\end{tabular}

\end{center}

Here $\tau_0^*=\frac{-b_0^*+\sqrt{-3}}{2}$, $\tau_i=\frac{-b_i+\sqrt{-3}}{2}$, where we choose $b_i, b_0, b_0^*$ such that $b_i\equiv b_0 \equiv b_0^* \equiv 1(6D)$, $b_i\equiv b_0 \equiv b_0^*(8)$, differing mod $\alpha$ as follows:

\begin{itemize}
	\item $b_0 \equiv\sqrt{-3}(\alpha)$	
	\item $b_0^*\equiv -\sqrt{-3}(\alpha)$
	\item $b_i-b_0\equiv i (m)$
\end{itemize}
We will also fix in all cases only $b_i$ such that $b_i\equiv \pm 1(8)$. This does not change the result, but simplifies some of the calculations. 
	
	We show in Lemma \ref{shim_1} in Section \ref{comp_shimura} that these traces are well-defined over the precise ring class fields we take, thus $X^*, U, T\in K$. Moreover, as each ratio is an algebraic integer from CM theory, we get further $X^*, U, T\in \OO_K$.

We note that, for $\alpha \equiv \pm \sqrt{-3}(4)$, the traces $U$ and $T$ also depend on $b_i (8)$, but we omit this from the notation. 
	
We define $u_{\alpha, b}=\begin{cases} 1&   \text{ if } \alpha\equiv 1(4)\\-1&   \text{ if } \alpha\equiv -1(4) \\ [\alpha]\left(\frac{-b}{4}\right)&  \text{ if } \alpha\equiv \pm\sqrt{-3}(4) \end{cases}$ and note that, for $m\equiv 3(4)$, the condition \eqref{cond} implies 
\begin{equation}\label{cond2}
u_{\alpha, -b_0^*}=-1.
\end{equation}

The goal of the section is to show:

\begin{Prop}\label{rel_3}
		\begin{enumerate}[(i)]
			\item $X^*=\frac{U}{L_{\overline{\alpha}}(1, \chi)}$, 
			\item $X^*\frac{ L_{\overline{\alpha}}(1, \chi)}{\alpha}= -u_{\alpha, -b_0^*}\eps(\sqrt{-3})\frac{T}{D'(1-\chi_D(\overline{\alpha})\eps(\overline{\alpha})\varphi(\overline{\alpha}))}$
		\end{enumerate}
	
\end{Prop}

	We also define the ideals prime to $m$:
\begin{itemize}	
	
	\item  $\A_j=\left(\tau_j\right)=[a_j, \tau_j]$, $j\not\equiv 0, 2b_0 (m)$, of norm $a_j=\frac{b^2_j+3}{4}$

	\item $\C_0'=(\tau_0)/(\alpha)=[c'_0, \tau_0]$ of norm $c_0'=\frac{{b_0}^2+3}{4m}$
	
	\item $\C_0^*=(\tau_0^*)/(\overline{\alpha})=[c_0^*, \tau_0^*]$ of norm $c_0^*=\frac{{b_0^*}^2+3}{4m}$
\end{itemize}

The idea of the proof is to write several sums of ratios of theta functions using Lemma \ref{sum_M} from the Section \ref{theta_M}, that states:
\[
\sum_{0\leq j \leq d-1}\Theta_M\left(z+\frac{j}{d}\right)=d\Theta_M(dz),\ \ d|3m^*.
\] 
We show that the ratios of theta functions in the sums are Galois conjugate under the Galois action given by the ideals $\A_j, \C'_0, \C^*_0$ defined above (via the Artin map). This is shown using Lemma \ref{inv1} from the Appendix. Taking the traces in the sums, and using the characters computed in the Appendix Lemma \ref{char} and \ref{charD}, we get relations between $X^*$, $T$ and $U$.

We will proceed now to computing these relations in Lemma \ref{rel_1} and Lemma \ref{rel_2}. First we show:

\begin{Lem}\label{rel_1} $X^*=-u_{\alpha, -b_0^*}\eps(\sqrt{-3})\frac{\alpha}{D}T+\frac{\chi_D(\overline{\alpha})\eps(\overline{\alpha})}{\varphi(\alpha)}(m-1)U.$

\end{Lem}

{\bf Proof:} For $\alpha\equiv \pm 1(4)$, using Lemma \ref{sum_M} for $d=m$ and $z=\frac{b_0+\sqrt{-3}}{6mD'}$, we have $\sum\limits_{0 \leq j \leq m-1} \Theta_M(\frac{-\overline{\tau_j}}{3mD'})= m \Theta_M(\frac{-\overline{\tau_j}}{3D'}).$ We can rewrite this as:
\[
\sum_{j} \frac{\Theta_M(-\overline{\tau_j}/3mD')}{\Theta_K(-\overline{\tau_j}/3)}
+
\frac{\Theta_M(-\overline{\tau_0}/3mD')}{\Theta_K(-\overline{\tau_0}/3)}
+
\frac{\Theta_M(-\overline{\tau^*_0}/3mD')}{\Theta_K(-\overline{\tau_0^*}/3)}
=
m \frac{\Theta_M(-\overline{\tau_j}/3D')}{\Theta_K(-\overline{\tau_j}/3)},
\]
where the sum is taken over all $j$ such that $b_j^2\not\equiv -3(m)$. Using the Appendix Lemma \ref{inv1} (iii), (iv), (v) we get on the LHS:
\[
\sum_j D'\sqrt{m}\left(\frac{\Theta_M(D'\tau_j)}{\Theta_K(\tau_j)}\right)^{\sigma_{\A_j}^{-1}} 
+
\frac{D'\sqrt{m}}{\varphi(\overline{\alpha})}\left(\frac{\Theta_M(D'\tau_0/m)}{\Theta_K(\omega)}\right)^{\sigma_{\C'_0}^{-1}}
+
\frac{D'\sqrt{m}}{\varphi(\alpha)}\left(\frac{\Theta_M(D'\tau_0^*/m)}{\Theta_K(\omega)}\right)^{\sigma_{\C^*_0}^{-1}}
\]
and multiplying by $\overline{\alpha}^{1/2}D^{1/3}$ we get:
\begin{small}
\[
\sum_j D'\overline{\alpha}\left(\frac{\Theta_M(D'\omega)}{\Theta_K(\omega)}D^{1/3}\alpha^{1/2}\right)^{\sigma_{\A_j}^{-1}} \eps(\A_j)\chi_D(\overline{\A_j})
+
\frac{D'\overline{\alpha}}{\varphi(\overline{\alpha})}\left(\frac{\Theta_M(D'\tau_0/m)}{\Theta_K(\omega)}\alpha^{1/2}D^{1/3}\right)^{\sigma_{\C'_0}^{-1}}\chi_D(\overline{\C'_0})\eps(\C'_0)
\]
\[
+
\frac{D'\overline{\alpha}}{\varphi(\alpha)}\left(\frac{\Theta_M(D'\tau_0^*/m)}{\Theta_K(\omega)}\alpha^{1/2}D^{1/3}\right)^{\sigma_{\C_0^*}^{-1}}\chi_D(\overline{\C^*_0})\eps(\C^*_0)
=
m\frac{\Theta_M(-\overline{\tau_j}/3D')}{\Theta_M(-\overline{\tau_j}/3)}D^{1/3}\overline{\alpha}^{1/2}.
\]
\end{small}
Taking the traces to $H_{3D'm}$, we get:
\[
D'\overline{\alpha}X\sum_j \eps(\A_j)\chi_D(\overline{\A_j})
+
\frac{D\overline{\alpha}}{\varphi(\overline{\alpha})} \chi_D(\overline{\C^*_0})\eps(\C^*_0)(m-1)U
=mT.
\]
Here we have used $\Tr_{H_{3D'm}/K}\left(\frac{\Theta_M(D'\tau_0/m)}{\Theta_K(\omega)}\alpha^{1/2}D^{1/3}\right)=0$ from Lemma \ref{shim_1}. We also note that the trace of $U$ is defined on $H_{3D'}/K$, hence the extra $(m-1)$-factor.

 We similarly obtain for $\alpha\equiv \pm \sqrt{-3}(4)$:
\[
X^{(i)}D\overline{\alpha}\sum_j \eps(\A_j)\chi_D(\overline{\A_j})
+
\frac{D\overline{\alpha}}{\varphi(\overline{\alpha})} \chi_D(\overline{\C^*_0})\eps(\C^*_0)(m-1)U
=mT.
\]

We computed the $\eps$ and $\chi_D$ characters in the Appendix. From Lemma \ref{charD}, $\chi_D(\overline{\A_j})=1$, $\chi_D(\overline{\C^*_0})=\chi_D(\overline{\alpha})$, and, from Lemma \ref{char}, $\eps(\C^*_0)=u_{\alpha, -b_0^*}\eps(\overline{\alpha})\eps(\sqrt{-3})$, thus we have 
\[
\sum_{i} \overline{\alpha}\eps(\A_j) D'X^* +\eps(\sqrt{-3})\frac{D'\overline{\alpha}\chi_D(\overline{\alpha})\eps(\overline{\alpha})}{\varphi(\alpha)}(m-1)U=mT.
\] 
As, from Lemma \ref{char}, $\eps(\A_j)=\left(\frac{(-b_j+b_0)/2}{m}\right)u'$, where $u'=1$ for $m\equiv 1(4)$ and $u'=\left[\frac{b_j}{4}\right]$ for $m\equiv 3(4)$,
we have $u'\sum_{j} \left(\frac{(-b_j+b_0)/2}{m}\right)=-u'\left(\frac{b_0}{m}\right)=-u'\left[\frac{\sqrt{-3}}{\alpha}\right]=
u'\eps(\sqrt{-3}) [\alpha][-\sqrt{-3}\alpha]$. This equals $-u_{\alpha, -b_j}\eps(\sqrt{-3})$.

Thus the final sum equals $ -\overline{\alpha}u_{\alpha, -b_j}\eps(\sqrt{-3})D'X +u_{\alpha, -b_j}\eps(\sqrt{-3})\frac{D'\overline{\alpha}\chi_D(\overline{\alpha})\eps(\overline{\alpha})}{\varphi(\alpha)}(m-1)U=mT,$ or equivalently $X=-u_{\alpha, -b_j}\eps(\sqrt{-3})\frac{\alpha}{D'}T+\frac{\chi_D(\overline{\alpha})\eps(\overline{\alpha})}{\varphi(\alpha)}(m-1)U$ for $\alpha\equiv \pm 1 (4)$.

For $\alpha \equiv \pm\sqrt{-3} (4)$, we get $-\overline{\alpha}u_{\alpha, -b_0^*}\eps(\sqrt{-3})DX^{(i)} +u_{\alpha, -b_0^*}\eps(\sqrt{-3})\frac{D\overline{\alpha}\chi_D(\overline{\alpha})\eps(\overline{\alpha})}{\varphi(\alpha)}(m-1)U=mT$.

\bigskip

We show a second linear relation between $X^*$, $U$ and $T$ below:

\begin{Lem}\label{rel_2}
$-u_{\alpha, -b_0^*}\eps(\sqrt{-3})\frac{\alpha}{D}T
+(m-1)U=mX^*$

\end{Lem}

{\bf Proof:} For $\alpha\equiv \pm1(4)$, using Lemma \ref{sum_M} for $d=m$ and $z=D'\frac{-b_0+\sqrt{-3}}{2m}$, we have $\sum\limits_{0 \leq j \leq m-1} \Theta_M(D'\frac{-b_j+\sqrt{-3}}{2m})= m \Theta_M(D'\omega).$ We rewrite this as:
\begin{small}
\[
\sum_{b_j} \frac{\Theta_M(D'\tau_j/m)}{\Theta_K(\omega)}\alpha^{1/2}D^{1/3}
+
\frac{\Theta_M(D'\tau_0/m)}{\Theta_K(\omega)}\alpha^{1/2}D^{1/3}
+
\frac{\Theta_M(D'\tau_0^*/m)}{\Theta_K(\omega)}\alpha^{1/2}D^{1/3}
= 
m 
\frac{\Theta_M(D'\omega)}{\Theta(\omega)}\alpha^{1/2}D^{1/3},
\]
\end{small} 

\noindent where the sum is taken over all $j$ such that $b_j^2\not\equiv -3(m)$. As $\frac{\Theta_M(D'\tau_i/m)}{\Theta_K(\omega)}
	=
	\frac{\sqrt{m}}{D'}\left(\frac{\Theta_M(-\overline{\tau_j}/3D')}{\Theta_K(-\overline{\tau_j}/3)}\right)^{\sigma_{\overline{\A_j}}^{-1}}$, from Lemma \ref{inv1}(vi), for the ideal $\overline{\A_j}=\left(-\overline{\tau_j}\right)$ of norm $a_j$, then the first sum equals 
	\[
	\frac{\alpha}{D'}\sum_{j} \left(\frac{\Theta_M(-\overline{\tau_j}/3D')}{\Theta_K(-\overline{\tau_j}/3)}\overline{\alpha}^{1/2}D^{1/3}\right)^{\sigma_{\overline{\A_j}}^{-1}} \chi_D(\A_j)\left[\frac{\overline{\alpha}}{\overline{\A_j}}\right].
	\]
From Lemma \ref{char} we have $\chi_D(\A_j)=1$, thus taking the traces, we get
\[
\frac{\alpha}{D'}T \sum_{j}\eps(\A_j)
+(m-1)U=mX,
\] 
where we have used Lemma \ref{shim_1} to show that $\Tr_{H_{3D'm}/K}\frac{\Theta_M(D'\tau_0/m)}{\Theta_K(\omega)}\alpha^{1/2}D^{1/3}=0$. We show similarly that $\frac{\alpha}{D'}T \sum_{j}\eps(\A_j)
+(m-1)U=mX^{(i)}$ for $\alpha \equiv \pm \sqrt{-3}(4)$.

We already computed in the previous lemma $\sum\limits_{j}\eps(\A_j)= -u_{\alpha, -b_0^*}\eps(\sqrt{-3})$. Thus we get $-u_{\alpha, -b_0^*}\eps(\sqrt{-3})\frac{\alpha}{D}T
+(m-1)U=mX^*$.

\bigskip

{\bf Proof Proposition \ref{rel_3}} From Lemma \ref{rel_1} and Lemma \ref{rel_2}, we get immediately Proposition \ref{rel_3} by solving the system
$\begin{cases}
X^*=-u_{\alpha, -b_0^*}\eps(\sqrt{-3})\frac{\alpha}{D'}T+\frac{\chi_D(\overline{\alpha})\eps(\overline{\alpha})}{\varphi(\alpha)}(m-1)U \\ 
mX^*=-u_{\alpha, -b_0^*}\eps(\sqrt{-3})\frac{\alpha}{D'}T
+(m-1)U\\
\end{cases}$. Thus $X^*=U(1-\frac{\chi_D(\overline{\alpha})\eps(\overline{\alpha})}{\varphi(\alpha)})$, or equivalently $L_{\alpha}(1, \chi)X^*=U$.
We also obtain 
$X^*(1-m\frac{\chi_D(\overline{\alpha})\eps(\overline{\alpha})}{\varphi(\alpha)})
=
-(1-\frac{\chi_D(\overline{\alpha})\eps(\overline{\alpha})}{\varphi(\alpha)})u_{\alpha, -b_0^*}\eps(\sqrt{-3})\frac{\alpha}{D}T$, which is equivalent to $X^*\frac{ L(1, \chi)}{\alpha}= -u_{\alpha, -b_0^*}\eps(\sqrt{-3})\frac{T}{D'(1-\chi_D(\overline{\alpha})\eps(\overline{\alpha})\varphi(\overline{\alpha}))}.$

\bigskip

 \subsection{Complex conjugates}\label{XYZ}
 In this section we will relate $U/\alpha$ to its complex conjugate.

We choose $b_{0, i} \equiv -1 (2D)$ for $i=0, 1, 2$ such that

	\begin{itemize}
		\item  $b_{0, i} \equiv -i (3)$
		\item $b_{0, i}\equiv \sqrt{-3}(\alpha)$
		\end{itemize}

We choose $b_{0, 0}\equiv b_{0, 1}\equiv b_{0, 2}(8) \equiv -b_0^*$, which is congruent to $\pm 1(8)$ for simplification of calculations. We note that then we have
\begin{equation}\label{b_0}
b_0^*\equiv -b_{0, 1} (24D\alpha),
\end{equation}
where $b_0^*$ was defined in Section \ref{XT}. Moreover, as $b_0^*\equiv -b_{0, 1}(8)$, we also have $u_{\alpha, -b_{0, i}}=u_{\alpha,  b_0^*}$.
\bigskip
For $\overline{\tau}_{0, i}=\frac{-b_{0, i}+\sqrt{-3}}{2}$, thus $-\overline{\tau}_{0, i}=\frac{b_{0, i}+\sqrt{-3}}{2}$, we define:

\begin{center}
\begin{tabular}{|l |l |l |}
\hline
& $\alpha \equiv \pm 1(4)$ &$\alpha \equiv \pm \sqrt{-3}(4)$ \\ 
\hline
$Y_i$
 &
$\Tr_{H_{3D'}/K}
\frac{\Theta_M\left(\frac{-\overline{\tau}_{0, i}}{3mD'}\right)}{\Theta_K(\omega)} \alpha^{1/2}D^{1/3}$
&
$
\Tr_{H_{6D}/K}
\frac{\Theta_M\left(\frac{-\overline{\tau}_{0, i}}{6mD}\right)}{\Theta_K(\omega)} \alpha^{1/2}D^{1/3}$ \\ 
\hline
$Y$ & $\Tr_{H_{3D'}/K} \frac{\Theta_M\left(\frac{-\overline{\tau}_{0, i}}{mD'}\right)}{\Theta_K (\omega)} \alpha^{1/2}D^{1/3},$ 
& $\Tr_{H_{6D}/K} \frac{\Theta_M\left(\frac{-\overline{\tau}_{0, i}}{2Dm}\right)}{\Theta_K (\omega)} \alpha^{1/2}D^{1/3}$ 
\\
\hline
$W_i$ &
$\Tr_{H_{3D'}/K}
\frac{\Theta_M\left(\frac{-\overline{\tau}_{0, i}}{3mD'}\right)}{\Theta_K\left(-\overline{\tau}_{0, i}/3\right)} \alpha^{1/2}D^{1/3}$  
& $
\Tr_{H_{6D}/K}
\frac{\Theta_M\left(\frac{-\overline{\tau}_{0, i}}{6mD}\right)}{\Theta_K\left(-\overline{\tau}_{0, i}/3\right)} \alpha^{1/2}D^{1/3}$ 
\\
\hline
\end{tabular}
\end{center}
\bigskip

Here we define $W_i$ for $i=1,2$ and $Y_i$ for $i=0,1,2$. We note that the traces are well defined from Lemma \ref{shim_1} and that $W_i, Y_i, Y\in \OO_K$, as they are traces of algebraic integers (again from CM-theory applied to modular functions evaluated at CM points).

We define the ideals:
\begin{itemize}
	\item $\C_0=\left(\tau_{0, 0}\right)/(\alpha\sqrt{-3})$ of norm $c_0$
	\item $\C_i=\left(\tau_{0, i}\right)/(\alpha)$ of norm $c_i$
	\item $\A_0=t_{\A}a+3m^*s'\tau_{0, 0}$, with $D' |t_{\A}$, of norm $a_0$
	\end{itemize}

\bigskip

 As $U$ will equal $Y/D'$ up to a sixth root of unity (see Corollary \ref{UY}), the goal of the section is to show Proposition \ref{conjU} in which we relate $U/\alpha$ to its complex conjugate via Proposition \ref{conjY} which related $Y/\alpha$ to its complex conjugate.

From Lemma \ref{sum_M} for $d=3$, we have 
\begin{equation}\label{sum_Y}
3Y=Y_0+Y_1+Y_2.
\end{equation}
The goal is to obtain a relation between $Y$ and $\overline{Y}$.  The general idea is to relate the traces to their complex conjugates by looking at various Galois conjugates, using the Galois action corresponding to the ideals $\A_0, \C_i$ defined above. These Galois conjugates are computed explicitly in Lemma \ref{inv1} in the Appendix. The characters $\eps$ and $\chi_D$ are computed in the Lemmas \ref{char} and \ref{charD} in the Appendix.

In Lemma \ref{WY} we relate the values of $W_i$ to $\overline{Y}$, which also implies Corollary \ref{form_Y}, relating the values of $Y_i$ to $\overline{Y}$ for $i=1, 2$. The essential tool is Lemma \ref{lemma_A} proved in the Section \ref{comp_shimura}, which allows us to freely move the $D'$ between the numerator and denominator in the traces. More precisely, for $\tau=\frac{-b_0+\sqrt{-3}}{2\cdot 3^e \cdot m^{e'}}$ with $e, e'\in \{0,1\}$ and $b_0\equiv \sqrt{3}(\alpha)$, we have:
\begin{equation}
\frac{\Theta_M(D'\tau)}{\Theta_K(\tau)}
=
\frac{1}{D'}\left(\frac{D}{3m^*}\right)\left(\frac{t_{\A}}{m^*}\right)\left(\frac{\Theta_M(\tau/D')}{\Theta_K(\tau)}\right)^{\sigma_{\A_{0}}^{-1}}. 
\end{equation}

Via a similar approach, we relate $Y_0$ to its complex conjugate in Lemma \ref{Y_0} and show that in certain cases $Y_0=0$.

 Plugging back in \eqref{sum_Y} the values for $Y_i$, we get Proposition \ref{conjY}, in which we show that $Y/\alpha$ is real or purely imaginary up to a third root of unity. Finally this gives us Proposition \ref{conjU} that $U/\alpha$ is real or purely imaginary up to a third root of unity.

\bigskip

We also note that we can define $U^{\circ}=U/\alpha$, $Y^{\circ}=Y/\alpha$, $Y_i^{\circ}=Y_i/\alpha$ for $i=0, 1, 2$ and $W_1^{\circ}=W_1/\alpha$. In particular, using \eqref{lemma_phi} we have $\Theta_K(-\tau_{0, i}/m)=\alpha \Theta_K(\omega)$ and $\Theta_K(-\tau_{0, i}/(3m))=\alpha \Theta_K(-\tau_{0, i}/3)$. Then, using Lemma \ref{OK} from Section \ref{comp_shimura}, $U^{\circ}, Y^{\circ}, Y_i^{\circ}, W_1^{\circ}$ are in $\OO_K$, and our results state:
\begin{itemize}
	\item $\overline{U}^{\circ}\stackrel{\zeta_6}{=}
U^{\circ}$, $U^{\circ}\in \OO_K$
	\item $\overline{Y}^{\circ}\stackrel{\zeta_6}{=}Y^{\circ}$, $Y^{\circ}\in \OO_K$
	\item $\overline{W_1}^{\circ}\stackrel{\zeta_6}{=}W_1^{\circ}$, $W_1^{\circ}\in \OO_K$,
\end{itemize}
where by $\stackrel{\zeta_6}{=}$ we mean equality up to a $6^{th}$ root of unity.
To ease the calculations, we will keep the notation $U, Y_i, W_i, Y$ in the proofs and statements below. We will denote the constant 
\[
c_{\alpha}=u_{\alpha, -b_{0}^*}\left(\frac{D}{3m^*}\right)\eps(\sqrt{-3})\eps(\overline{\alpha})\chi_D(\overline{\alpha})\frac{\varphi(\alpha)}{\alpha},
\]
which equals a sixth root of unity. We also note that, under \eqref{cond}, from Lemma \ref{eps_1} we have $u_{\alpha, -b_0^*}\eps(\sqrt{-3})\eps(\overline{\alpha})=1$. Thus actually:
\[
c_{\alpha}=\left(\frac{D}{3m^*}\right)\chi_D(\overline{\alpha})\frac{\varphi(\alpha)}{\alpha}.
\]

We also recall that we computed the value of the character $\chi_{\omega}(D)$ in Section \ref{cubic}.

\bigskip

We start by relating $W_1$ to $\overline{U}$:

\begin{Lem}\label{WU} $\ds \frac{W_1}{\alpha}=D'\left(\frac{D}{3m^*}\right)c_{\alpha}\frac{\overline{U}}{\overline{\alpha}}$.
\end{Lem}
{\bf Proof:} From Lemma \ref{inv1}(ii) from the Appendix, we get $\ds W_1=D'\frac{\alpha}{\varphi(\overline{\alpha})}\chi_D(\overline{\C_1})\eps(\overline{\C_1})\overline{U}$. Then, using the characters computed in Lemma \ref{char} and Lemma \ref{charD} in the Appendix, we further compute $W_1=D'\frac{\alpha}{\varphi(\overline{\alpha})}\chi_D(\overline{\alpha})u_{\alpha, -b_0^*}\eps(\sqrt{-3})\eps(\overline{\alpha})\overline{U}$.

\bigskip

We will now relate $W_i$ to $\overline{Y}$ below:

\begin{Lem}\label{WY} 
$\ds \frac{W_1}{\alpha}= \chi_{\omega}(D)c_{\alpha} \frac{\overline{Y}}{\overline{\alpha}}$ 
and 
$\ds\frac{W_2}{\alpha}= c_{\alpha}\frac{\overline{Y}}{\overline{\alpha}}$

\end{Lem}

{\bf Proof:} For $\alpha\equiv \pm 1(4)$,  from Lemma \ref{inv1} (ii) in the Appendix, we have $\ds \frac{\Theta_M\left(-\overline{\tau_{0, i}}/3mD'\right)}{\Theta_K\left(-\overline{\tau_{0, i}}/3\right)}
=
\frac{D'\sqrt{m}}{\varphi(\overline{\alpha})}\left(\frac{\Theta_M\left(D'\tau_{0, i}/m\right)}{\Theta_K(\omega)}\right)^{\sigma_{\C_i}^{-1}}$. We use Lemma \ref{lemma_A} to rewrite $\frac{\Theta_M\left(D'\tau_{0, i}/m\right)}{\Theta_K(\omega)}
= 
\frac{1}{D'}\left(\frac{t_{\A}}{m^*}\right)\left(\frac{D}{3m^*}\right)\left(\frac{\Theta_M\left(\tau_{0, i}/mD'\right)}{\Theta_K(\omega)}\right)^{\sigma_{\A_0}^{-1}}$ 
for the ideal $\A_0=(t_{\A}a+3s'm^*\frac{-b_{0}+\sqrt{-3}}{2})$, $D'|t_{\A}$. Thus we get:
\[
\frac{\Theta_M\left(-\overline{\tau_{0, i}}/3mD'\right)}{\Theta_K\left(-\overline{\tau_{0, i}}/3\right)}\alpha^{1/2}D^{1/3}
=
\left(\frac{t_{\A}}{m^*}\right)\left(\frac{D}{3m^*}\right)\frac{\alpha}{\varphi(\overline{\alpha})}\left(\frac{\Theta_M\left(\tau_{0, i}/mD'\right)}{\Theta_K(\omega)} \overline{\alpha}^{1/2} D^{1/3}\right)^{\sigma_{\A_0\C_i}^{-1}} \eps(\overline{\A_0\C_i}) \chi_D(\overline{\A_0\C_i}).
\]
Taking the traces, this is $
W_i = u \overline{Y}$, where
\[
u
=
\left(\frac{t_{\A}}{m^*}\right)\left(\frac{D}{3m^*}\right)\frac{\alpha}{\varphi(\overline{\alpha})}\eps(\overline{\A_0\C_i}) \chi_D(\overline{\A_0\C_i}).
\]
Similarly we obtain $W_i = u \overline{Y}$ for $\alpha \equiv \pm \sqrt{-3}(4)$.

Then from Lemma \ref{char} and \ref{charD}, where we computed the characters, we get:
\[
\eps(\overline{\A_0\C_i}) \chi_D(\overline{\A_0\C_i})
=u_{\alpha, b_{0, i}}\left(\frac{t}{m}\right)\eps(\overline{\alpha})\chi_{\omega}(D)^{2i-1}\chi_D(\overline{\alpha})
\]
Thus $W_i= u_{\alpha, b_{0, i}}\chi_{\omega}(D)^{2i-1}\eps(\sqrt{-3})\left(\frac{D}{3m^*}\right)\frac{\eps(\overline{\alpha})\chi_D(\overline{\alpha}) \alpha}{\varphi(\overline{\alpha})} \overline{Y},$ which gives us the result of the lemma.
	
\bigskip

From the above two lemmas, we get immediately: 

\begin{cor}\label{UY} $\ds D'U=\chi_{\omega^2}(D) \left(\frac{D}{3m^*}\right)Y$

\end{cor}

Our goal is now to relate $Y/\alpha$ to its complex conjugate. We first note:
\begin{Lem}\label{sum_3}
$3Y=Y_0+Y_1+Y_2.$
\end{Lem}

{\bf Proof:} For $\alpha\equiv \pm 1(4)$, we apply Lemma \ref{sum_M} for $d=3$ and $-\overline{\tau_{0, 1}}/3D'm$ and get:
\[
\Theta_M\left(-\overline{\tau_{0, 0}}/3D'm\right)
+
\Theta_M\left(-\overline{\tau_{0, 1}}/3D'm\right)
+
\Theta_M\left(-\overline{\tau_{0, 2}}/3D'm\right)
=
3\Theta_M\left(-\overline{\tau_{0, 1}}/D'm\right).
\]

Multiplying by $\frac{D^{1/3}\alpha^{1/2}}{\Theta(\omega)}$ and taking traces we get $3Y=Y_0+Y_1+Y_2.$ Similarly we show the same relation for $\alpha \equiv \pm \sqrt{-3}(4)$.

\bigskip

The goal is now to relate each of the terms $Y_i$ either to $\overline{Y}$, or to their own complex conjugate. Noting that $Y_1=(\omega^2 \sqrt{-3})W_1$ and $Y_2=(-\omega \sqrt{-3})W_2$, we get immediately from Lemma \ref{WY}:

\begin{cor}\label{form_Y} 
\begin{enumerate}[(i)]
	\item $\ds\frac{Y_1}{\alpha}=(\omega^2\sqrt{-3})\chi_{\omega}(D)c_{\alpha}\frac{\overline{Y}}{\overline{\alpha}},$
	
	\item	
	$\ds\frac{Y_2}{\alpha}=(-\omega\sqrt{-3})c_{\alpha} \frac{\overline{Y}}{\overline{\alpha}}$
	
	\end{enumerate}
\end{cor}


Using similar methods as in Lemma \ref{WY},  we relate below $Y_0$ to its complex conjugate:

\begin{Lem}\label{Y_0}
\begin{enumerate}[(i)]
	\item  $Y_0 = 0$, if $D\equiv \pm 1, \pm 4 (9)$,
	\item  If $D\equiv \pm 2$, let $\zeta=\left(\frac{b+\sqrt{-3}}{2}\right)/(\sqrt{-3})$, $b\equiv 0(9)$, where $b\equiv -1(2D)$, $b\equiv \sqrt{-3}(\alpha)$, we have
\[
\frac{Y_0}{\alpha}
=
-\omega^2\chi_D(\zeta)\eps(\sqrt{-3})c_{\alpha}
\frac{\overline{Y}_0}{\overline{\alpha}},
\]

\end{enumerate}

\end{Lem}

{\bf Proof:} For $\alpha\equiv \pm 1(4)$, from Lemma \ref{inv1} (i) in the Appendix, we have $\frac{\Theta_M\left(-\overline{\tau_{0, 0}}/3D'm\right)}{\Theta_K\left(-\overline{\tau_{0, 0}}\right)}
=
\frac{\varphi(\alpha) D'}{3\sqrt{m}}\left(\frac{\Theta_M\left(D'\tau_{0, 0}/3m\right)}{\Theta_K\left(\tau_{0, 0}/9\right)}\right)^{\sigma_{\C_0}^{-1}}.$ Then for the ideal $\A_0=t_{\A}a+3m^*s'\tau_{0, 0}$, with $D'|t_{\A}$, from Lemma \ref{lemma_A} we get $\frac{\Theta_M\left(D'\tau_{0, 0}/3m\right)}{\Theta_K\left(\omega\right)}
=
\frac{\left(\frac{D}{3m^*}\right)}{D'}\left(\frac{t_{\A}}{m^*}\right)\left(\frac{\Theta_M\left(\tau_{0, 0}/3mD'\right)}{\Theta_K\left(\omega\right)}\right)^{\sigma_{\A_0}^{-1}},$ which gives us:
\[
\frac{\Theta_M\left(-\overline{\tau_{0, 0}}/3D'm\right)}{\Theta_K\left(\omega\right)}\alpha^{1/2}D^{1/3}
=
\frac{3\varphi(\alpha)}{\overline{\alpha}}\left(\frac{D}{3m^*}\right)\left(\frac{t_{\A}}{m^*}\right) \frac{\Theta_K\left(\omega\right)}{\Theta_K\left(\tau_{0, 0}/9\right)}\left(\frac{\Theta_M\left(\tau_{0, 0}/3mD'\right)}{\Theta_K\left(\omega\right)}\right)^{\sigma_{C_0}^{-1}\sigma_{\A_0}^{-1}}\overline{\alpha}^{1/2}D^{1/3}.
\]
Taking the trace from $H_{3D'}$ to $K$ we get $Y_0=
u\overline{Y_0}$,
where 
\[
u= \frac{3\varphi(\alpha)}{\overline{\alpha}}\left(\frac{D}{3m^*}\right)\left(\frac{t_{\A}}{m^*}\right) \frac{\Theta_K\left(\omega\right)}{\Theta_K\left(\tau_{0, 0}/9\right)} \chi_D(\overline{\A_0\C_0})\eps(\overline{\A_0\C_0}).
\] Similarly, we compute 
$Y_0=u\overline{Y_0}$ for $\alpha\equiv \pm \sqrt{-3} (4)$ using again Lemma \ref{inv1} (i) and Lemma \ref{lemma_A}.

We compute now $u$. In Lemma \ref{char} and \ref{charD} we computed the characters $\eps$ for $\A_0$ and $\C_i$. The character $\chi_D(\C_0)$ is the only one that depends on $b_0 \mod 9$. Then, if $Y_0\neq 0$, we must get a constant for $\Theta_K(\tau_{0}/9)\chi_D((\frac{b_{0, 0}/3+\sqrt{-3}}{2}))$ when we vary $b_{0, 0}\equiv 0, 3, -3 \mod 9$. Let $b_0\equiv 3, b_0'\equiv -3, b_0''\equiv 0 (9)$, $b_0 \equiv b_0'\equiv b_0'' \equiv -1(2D).$

We note that in all cases we have $\frac{-b_{0, 0}/3\sqrt{-3}+1}{2}\equiv -\frac{\sqrt{-3}}{3} \omega (D)$, and $\frac{-b_{0}/3\sqrt{-3}+1}{2}\equiv -\omega (3)$,  $\frac{-b'_{0}/3\sqrt{-3}+1}{2}\equiv -\omega^2 (3)$ and $\frac{-b''_{0}/3\sqrt{-3}+1}{2}\equiv -1 (3)$. Denote $\zeta=\frac{-b''_{0}/3\sqrt{-3}+1}{2}$. Then
\begin{itemize}
	\item  $\chi_D(\frac{-b_{0}/3\sqrt{-3}+1}{2})=\chi_{\omega^2}(D)\chi_D(\zeta)$,
	\item $\chi_D(\frac{-b'_{0}/3\sqrt{-3}+1}{2})=\chi_{\omega}(D)\chi_D(\zeta)$,
	\item $\chi_D(\frac{-b''_{0}/3\sqrt{-3}+1}{2})=\chi_D(\zeta)$.
\end{itemize}	
	
	As $\Theta_K((-3+\sqrt{-3})/18)=-3\omega \Theta(\omega)$ and $\Theta_K((9+\sqrt{-3})/18)=-3\Theta(\omega)$, we must have $\omega^2\chi_D((\frac{-b_{0}/3\sqrt{-3}+1}{2}))
=
\omega\chi_D((\frac{-b'_{0}/3\sqrt{-3}+1}{2}))
=
\chi_D((\frac{-b''_{0}/3\sqrt{-3}+1}{2}))
$. Thus $Y_0\neq 0$ only if $\chi_{\omega}(D)=\omega^2$. Using the values of Lemma \ref{char} and \ref{charD}:
\[
\chi_D(\overline{\A_0})\eps(\overline{\A_0\C_0})=\chi_{\omega^2}(D)\left(\frac{t}{m}\right) \eps(\overline{\alpha})
u_{\alpha, b_{0, 0}}
\]
Thus $Y_0
=
-u_{\alpha, b_{0, 0}}\left(\frac{D}{3m^*}\right) \frac{\varphi(\alpha)}{\overline{\alpha}}
\chi_D(\alpha)\eps(\overline{\alpha})\chi_{\omega}(D)\chi_D(\zeta)
\overline{Y}_0$ for $\chi_{\omega}(D)=\omega^2$ and $Y_0=0$ in the remaining cases.

\bigskip

\bigskip

Using Lemma \ref{sum_3}, Lemma \ref{Y_0} and Corollary \ref{form_Y}, we can finally relate $Y/\alpha$ to its complex conjugate:
\begin{Prop}\label{conjY} $\ds \frac{\overline{Y}}{\overline{\alpha}}= c_{\alpha} d_{\alpha}\frac{Y}{\alpha}$, where $d_{\alpha}= 
\begin{cases}
1 & \text{ if } D\equiv \pm 1(9) \\
-\omega^2 & \text{ if } D\equiv \pm 4(9) \\
-\omega^2\eps(\sqrt{-3})\chi_D(\overline{\zeta}) & \text{ if } D\equiv \pm 2(9),\\
\end{cases}.$

\end{Prop}

{\bf Proof:}  From Corollary \ref{form_Y}, we compute $Y_1+Y_2=(\omega^2\chi_{\omega}(D)-\omega)\sqrt{-3} c_{\alpha}\frac{\alpha}{\overline{\alpha}}\overline{Y}.$
Thus we have in Lemma \ref{sum_3}:
\begin{equation}
3Y=Y_0+cc_\alpha\frac{\alpha}{\overline{\alpha}}\overline{Y},  
\end{equation}
where $c=(\omega^2\chi_{\omega}(D)-\omega)\sqrt{-3}$.

For $\chi_{\omega}(D)=1, \omega$, we have from Lemma \ref{Y_0} that $Y_0=0$ in these cases, hence:

\begin{enumerate}[(i)]
	\item $\chi_{\omega}(D)=1$: $c=3$, thus $\ds \frac{Y}{\alpha}=c_{\alpha}\frac{\overline{Y}}{\overline{\alpha}}$.		
	\item  $\chi_{\omega}(D)=\omega$: $c=-3\omega^2$, thus $\ds \frac{Y}{\alpha}=-\omega^2c_{\alpha}\frac{\overline{Y}}{\overline{\alpha}}$.
\end{enumerate}

Finally for $\chi_{\omega}(D)=\omega^2$ we have $c=0$ and thus $3Y=Y_0$. From Lemma \ref{Y_0},  $Y_0
=
-\omega^2 \chi_D(\zeta)\eps(\sqrt{-3})c_{\alpha}\frac{\alpha}{\overline{\alpha}} \overline{Y}_0$, thus we get $\frac{Y}{\alpha}
=
-\omega^2\chi_D(\zeta)\eps(\sqrt{-3})\frac{\overline{Y}}{\overline{\alpha}}. $

\bigskip	

We could directly compute now how $\overline{U}/\overline{\alpha}$ relates to $U/
\alpha$ using Lemma \ref{UY}. However, for the sake of completeness, we relate $W_1$ to $U$ directly and compute the complex conjugate of $W_1/\alpha$. We do this below:

\begin{cor}\label{conjW} $\ds \frac{\overline{W_1}}{\overline{\alpha}}= c_{\alpha}t_{\alpha} \frac{W_1}{\alpha}$
for $t_{\alpha}= 
\begin{cases}
1 & \text{ if } D\equiv \pm 1(9) \\
-1 & \text{ if } D\equiv \pm 4(9) \\
-\eps(\sqrt{-3})\chi_D(\overline{\zeta}) & \text{ if } D\equiv \pm 2(9),\\
\end{cases}.$

\end{cor}
	
{\bf Proof:} From Lemma \ref{WY} and Proposition \ref{conjY}, we compute
	
\begin{enumerate}[(i)]
	\item $\chi_{\omega}(D)=1$: $W_1=Y$, thus $W_1= u_{\alpha, b_{0, 1}}c'\frac{\eps(\overline{\alpha})\chi_D(\overline{\alpha})\alpha}{\varphi(\overline{\alpha})}	\overline{W}_1.$

	\item  $\chi_{\omega}(D)=\omega$: $W_1=-\omega Y$, which implies $\overline{W}_1=-\omega^2\overline{Y}$ and thus $W_1=-u_{\alpha, b}\omega c'\frac{\eps(\overline{\alpha})\chi_D(\overline{\alpha})\alpha}{\varphi(\overline{\alpha})}
 \overline{W_1}.$

	\item $\chi_{\omega}(D)=\omega^2$: $W_1=-\omega^2\chi_D(\zeta)\eps(\sqrt{-3}) Y$ and thus $W_1=-u_{\alpha, b}\left(\frac{D}{3m^*}\right)\chi_D(\overline{\zeta}) \frac{\varphi(\alpha)}{\overline{\alpha}}
\chi_D(\overline{\alpha})\eps(\overline{\alpha})\overline{W_1}.$

\end{enumerate}

\bigskip

Combining this with Lemma \ref{WU}, we get the relation between $U$ and $W_1$:
\begin{cor}\label{WU2}
$W_1=t_{\alpha}\left(\frac{D}{3m^*}\right)D'U$.
\end{cor}

{\bf Proof:} We recall from Lemma \ref{WU} that $W_1=D'\frac{\alpha}{\varphi(\overline{\alpha})}\chi_D(\overline{\alpha})u_{\alpha, -b_0^*}\eps(\sqrt{-3})\eps(\overline{\alpha})\overline{U}$.  Combining with Proposition \ref{conjY}, we get $W_1=t_{\alpha}\left(\frac{D}{3m^*}\right)D'U$.
\bigskip

Finally we relate the complex conjugate of $U/\alpha$ to itself below:  

\begin{Prop}\label{conjU}
We have $\ds \frac{U}{\alpha}=c_{\alpha}t'_{\alpha}\frac{\overline{U}}
{\overline{\alpha}}$, for  $t'_{\alpha}= t_{\alpha}^3=
\begin{cases}
1 & \text{ for } D\equiv \pm 1(9) \\
-1 & \text{ for } D\equiv \pm 4(9) \\
-\eps(\sqrt{-3})& \text{ for } D\equiv \pm 2(9)\\
\end{cases}.
$
\end{Prop}

{\bf Proof:} From Corollary \ref{WU2}, we have $W_1=t_{\alpha}\left(\frac{D}{3m^*}\right)D'U$.   Thus applying Proposition \ref{conjW}, we get the result $\frac{U}{\alpha}=t^3_{\alpha} \left(\frac{D}{3m^*}\right)\eps(\overline{\alpha})\chi_D(\overline{\alpha}) \eps(\sqrt{-3})\frac{\varphi(\alpha)}{\alpha}\frac{\overline{U}}{\overline{\alpha}}$

\bigskip


\subsection{Proof of Theorem \ref{thm_square}}\label{proof} Finally, we are ready to prove Theorem \ref{thm_square}. 

We define:
\begin{equation}\label{def_Z}
Z=\chi_D(\overline{\alpha})\chi_3(\overline{\alpha})\Tr_{H_{3D^*}/K} 
 \left(\frac{\Theta_M(D^*\tau_0/m^*)}{\Theta_K(\omega)}D^{1/3} \alpha^{-1/2}\right),
\end{equation}
where $\chi_3(\alpha):=\frac{\alpha\eps(\sqrt{-3})}{\varphi(\alpha)}$ is the cubic root of unity $\omega^r$ such that $\alpha/\omega^r=\pm \varphi(\alpha)$.

From the property \eqref{lemma_phi}, $\Theta_K(\tau_0/m)=\varphi(\alpha)\Theta_K(\omega)$, as $(\overline{\alpha})=[m, \tau_0]$. Thus we can rewrite
\begin{equation}
Z=\chi_D(\overline{\alpha})\eps(\sqrt{-3})\chi_3(\alpha)\Tr_{H_{3D^*}/K} 
 \left(\frac{\Theta_M(D^*\tau_0/m^*)}{\Theta_K(\tau_0/m)}D^{1/3} \alpha^{1/2}\right),
\end{equation}
and we show in Proposition \ref{OK} in Section \ref{comp_shimura} that $Z\in \OO_K$. We have defined $Z$ such that 
\begin{equation}\label{ZU}
Z=\frac{U}{\alpha}\chi_D(\overline{\alpha})\chi_3(\overline{\alpha}),
\end{equation}
in order to have from Proposition \ref{conjU} 
\begin{equation}
\overline{Z}=t'_{\alpha}u_{\alpha, -b_0^*}\left(\frac{D}{3m^*}\right)\eps(\overline{\alpha})Z,
\end{equation}
giving us $\overline{Z}=\pm Z$.

It is immediate to see \eqref{ZU} for $m\equiv 1(4)$ from \eqref{def_Z}. For $m\equiv 3(4)$, we first note that $\frac{U}{\alpha}=Z'^{(i)}$, where $Z'^{(i)}=\Tr_{H_{6D}/K} \frac{\Theta_M\left(D\tau_i/(2m)\right)}{\Theta_K (\tau_i/m)} \alpha^{1/2}D^{1/3}$. In Section \ref{XZ_i} we will show in Lemma \ref{X_1} that $Z'=Z'^{(i)}$ under the condition \eqref{cond}, where $Z'=\Tr_{H_{12D}/K} \frac{\Theta_M\left(D\tau_i/m\right)}{\Theta_K (\tau_i/m)} \alpha^{1/2}D^{1/3}$. Thus $Z'=Z'^{(i)}=U/\alpha$ and we again have $Z=\frac{U}{\alpha}\chi_D(\overline{\alpha})\chi_3(\overline{\alpha}).$

Now we are ready to show:

\begin{cor}\label{real} For $Z$ defined above, we have:
\[
\frac{L(E_{D, \alpha}, 1)}{\Omega_{D\alpha}\overline{\Omega}_{D\alpha}}=c_0(2^eZ/\sqrt{-3})^2,
\]
with $\overline{Z}=-c_0Z$ and $c_0=(-1)^{v_3(c_{E_{D, \alpha}})}$, $e=1$ for $m\equiv 1(4)$ and $e=0$ for $m\equiv 3(4)$.


\end{cor}

{\bf Proof:} We recall, from Corollary \ref{rel_3}, $X^*
=
\frac{U}{L_{\overline{\alpha}}(1, \chi)}$ and from Theorem \ref{thm1} we have $\frac{L(E_{D, \alpha}, 1)}{\Omega_{D\alpha}\overline{\Omega}_{D\alpha}}=|2L_{\overline{\alpha}}(1, \chi)\frac{X^*}{\sqrt{-3}\alpha}|^2$, thus $\frac{L(E_{D, \alpha}, 1)}{\Omega_{D\alpha}\overline{\Omega}_{D\alpha}}=\frac{4}{3}|\frac{U}{\alpha}|^2$ and further $\frac{L(E_{D, \alpha}, 1)}{\Omega_{D\alpha}\overline{\Omega}_{D\alpha}}=|2Z|^2/3$.

Using Proposition \ref{conjU} we have $\ds \overline{Z}=t'_{\alpha}u_{\alpha, -b_0^*}\left(\frac{D}{3m^*}\right)\eps(\overline{\alpha})Z$, which is equivalent to:
\begin{itemize}
	\item for $\alpha \equiv 1(4)$:  $\overline{Z}=t_{\alpha}'\left(\frac{D}{3m^*}\right)\eps(\overline{\alpha})Z=t_{\alpha}'\left(\frac{D}{3m^*}\right)\eps(\sqrt{-3})Z$, where we have used Lemma \ref{eps_1}.

	\item for $\alpha \equiv -1(4)$: $\overline{Z}=-t_{\alpha}'\left(\frac{D}{3m^*}\right)\eps(\overline{\alpha})Z$. As $\eps(\overline{\alpha})=-\eps(\sqrt{-3})$ from Lemma \ref{eps_1}, we get $\overline{Z}=t_{\alpha}'\left(\frac{D}{3m^*}\right)\eps(\sqrt{-3})Z$
	
	\item for $\alpha \equiv \pm \sqrt{-3}(4)$: $\overline{Z}=-t_{\alpha}'\eps(\sqrt{-3})\left(\frac{D}{3m^*}\right)\eps(\overline{\alpha})Z$ for $\alpha \equiv \pm \sqrt{-3} (4)$, as $u_{\alpha, -b_0^*}=-1$ from \eqref{cond}. As $\eps(\sqrt{-3})=-\eps(\overline{\alpha})$ from Lemma \ref{eps_1}, we get $\overline{Z}=t_{\alpha}'\left(\frac{D}{3m^*}\right)\eps(\sqrt{-3})Z$.
	
\end{itemize}
Thus $\overline{Z}=-cZ$ for $c=\begin{cases} -\left(\frac{D}{3m^*}\right)\left[\frac{\alpha}{\sqrt{-3}}\right], & D\equiv \pm 1\\ \left(\frac{D}{3m^*}\right)\left[\frac{\alpha}{\sqrt{-3}}\right], & D\equiv \pm 4 (9) \\ \left(\frac{D}{3m^*}\right), & D\equiv \pm 2(9) \end{cases},$ and this is exactly $c_0$ (see Appendix \ref{tam}).

\bigskip

As $Z\in \OO_K$, we get immediately:

\begin{cor}\label{integral} $Z\in \ZZ$ if $c_0=-1$, and $Z/\sqrt{-3}\in\ZZ$ if $c_0=1$.
\end{cor}


\section{Shimura reciprocity computations}\label{comp_shimura}

\subsection{Fields of definition}

In the current section, we will show that all the traces used in Section \ref{X} are indeed well-defined. This is proved below by applying Shimura's reciprocity law.

\begin{Lem}\label{shim_1} For $\tau_0=\frac{-b_0+\sqrt{-3}}{2}$, $b_0\equiv \sqrt{-3}(\alpha)$, $e=\pm 1, e', e''\in\{0, 1\}$:
\begin{enumerate}[(i)]
	\item $\frac{\Theta_M(D'^e\frac{\tau_0}{3^{e'}m})}{\Theta_K(\tau_0/m^{e''})}\overline{\alpha}^{1/2}D^{1/3}\in H_{3D^*}$

	\item $\frac{\Theta_M(D'^e\frac{\tau_0}{3^{e'}})}{\Theta_K(\tau_0)}\overline{\alpha}^{1/2}D^{1/3}\in H_{3D'm^*}$

	\item $\Tr_{H_{3mD^*}/K} \frac{\Theta_M(D'^e\frac{\tau_0}{3^{e'}m})}{\Theta_K(\tau_0)}\alpha^{1/2}D^{1/3}=0$

\end{enumerate}
	
\end{Lem}

{\bf Proof:} Let $\ds G(z)=\frac{\Theta_M\left(D'^e\frac{z}{3^{e'}m}\right)}{\Theta_K(\frac{z}{m^{e''}})}$. For (i) it is enough to show the invariance under $\Gal(K^{ab}/H_{3D^*})$. Thus let $\A$ be an ideal of norm $a$ prime to $3D^*m$ with generator $k_{\A}=ta+s\tau_0$ such that $3D^*|s$. We want to show that $G(\tau_0)^{\sigma_{\A}^{-1}}=G(\tau_0)$. Using Lemma \ref{Shimura_A}, we have:
\[
G(\tau_0)^{\sigma_{\A}^{-1}}=
G(\left(\begin{smallmatrix} ta-sb & -scm\\  s& t\end{smallmatrix}\right)\tau_0),
\]
 where $c=\frac{b_0^2}{4ma}$. We compute $\Theta_M(\left(\begin{smallmatrix} D'^e & 0 \\ 0 & 3^{e'}m\end{smallmatrix}\right)\left(\begin{smallmatrix} ta-sb & -scm \\  s & t\end{smallmatrix}\right)\tau_0)=\Theta_M(\left(\begin{smallmatrix} ta-sb & -scD'^e/3^{e'}  \\  3^{e'}ms/D'^e & t\end{smallmatrix}\right)\frac{D'^e\tau_0}{3^{e'}m})$. Then $3m^*|3^{e'}ms/D'^e$ and we get 
 \[
 \Theta_M(\left(\begin{smallmatrix} D'^e & 0 \\ 0 & 3^{e'}m\end{smallmatrix}\right)\left(\begin{smallmatrix} ta-sb & -scm \\  s & t\end{smallmatrix}\right)\tau_0)=\left(\frac{t}{3m^*}\right)(s\tau_0+t)\Theta_M(\tau_0),
 \]
 where we have used Lemma \ref{modular} for the matrix $\left(\begin{smallmatrix} ta-sb & -scD'^e/3^{e'}  \\  3^{e'}ms/D'^e & t\end{smallmatrix}\right)$ in $\Gamma_0(3m^*)$.
 
 Similarly we get $\Theta_K(\left(\begin{smallmatrix} 1 & 0 \\ 0 & 3^{e'}m^{e''}\end{smallmatrix}\right)\left(\begin{smallmatrix} ta-sb & -scm \\  s & t\end{smallmatrix}\right)\tau_0)=\left(\frac{t}{3}\right)(s\tau_0+t)\Theta_K(\tau_0/3^{e'})$ using \eqref{modular2}. Thus we get:
 \[
 G(\tau_0)^{\sigma_{\A}^{-1}}=\left(\frac{t}{m^*}\right)G(\tau_0).
 \]
 
We compute $(\overline{\alpha}^{1/2})^{\sigma_{\A}^{-1}}=\left(\frac{ta-sb}{m^*}\right)$. For $\alpha\equiv 1(4)$ this is immediate from the reciprocity law $\left[\frac{\overline{\alpha}}{\A}\right]=\left[\frac{\tau_0}{\overline{\alpha}}\right]=\left[\frac{ta-sb}{\overline{\alpha}}\right]=\left(\frac{ta-sb}{m}\right)$. For $\alpha\equiv -1, \pm \sqrt{-3}(4)$, we have $ta+s\tau_0\equiv ta-sb \equiv \pm 1(4)$. We compute from \eqref{rec_law0} $\left[\frac{\overline\alpha}{\A}\right]=\left[\frac{ta-sb}{\overline{\alpha}}\right][\overline{\alpha}][(ta-sb) \overline{\alpha}][(ta-sb)]=\left(\frac{ta-sb}{m^*}\right)$. As $t(ta-sb) \equiv 1(a)$ and $t(ta-sb) \equiv 1(m^*/m)$, we get $G(\tau_0)\overline{\alpha}^{1/2}$ is invariant under $\A$, from which we get $(i)$. The proof of $(ii)$ is similar.

For (iii), we compute $(\alpha^{1/2})^{\sigma_{\A}^{-1}}=\left[\frac{\alpha}{\A}\right]\alpha^{1/2}$ and we get $\left[\frac{\alpha}{\A}\right]=\left(\frac{ta}{m^*}\right)$. Thus $(G(\tau_0)\alpha^{1/2})^{\sigma_{\A}^{-1}}=\left(\frac{a}{m}\right)G(\tau_0)\alpha^{1/2}$, which will have trace $0$.  

\bigskip

Similarly, one can show:

\begin{Lem}\label{shim_3} Under the same conditions as above, for $m\equiv 3 (4)$ we have:
\begin{enumerate}[(i)]
	\item $\frac{\Theta_M(D^e\frac{\tau_0/2}{3^{e'}m})}{\Theta_K(\tau_0/m^{e''})}\overline{\alpha}^{1/2}D^{1/3}\in H_{6D},$ $\frac{\Theta_M(D^e\frac{\tau_0}{3^{e'}m})}{\Theta_K(\tau_0/m^{e''})}\overline{\alpha}^{1/2}D^{1/3}\in H_{12D},$

	\item $\frac{\Theta_M(D^e\frac{\tau_0/2}{3^{e'}})}{\Theta_K(\tau_0)}\overline{\alpha}^{1/2}D^{1/3}\in H_{6Dm},$ $\frac{\Theta_M(D^e\frac{\tau_0}{3^{e'}})}{\Theta_K(\tau_0)}\overline{\alpha}^{1/2}D^{1/3}\in H_{12Dm},$	
		
	\item $\Tr_{H_{6mD}/K} \frac{\Theta_M(D^e\frac{\tau_0/2}{3^{e'}m})}{\Theta_K(\tau_0)}\alpha^{1/2}D^{1/3}=0$

\end{enumerate}
\end{Lem}

Thus we get $X^*$ and $T$ well-defined from (ii) of Lemma \ref{shim_1} and Lemma \ref{shim_3}. From (i) of Lemma \ref{shim_1} and Lemma \ref{shim_3} for $e''=0$, $U, Y_i, Y, W_i$ are well-defined, while for $e"=1$ then $Z, U^{\circ}, Y_i^{\circ}, Y^{\circ}, W_i^{\circ}$ are well-defined.
\bigskip

We finally remark that $F(\tau_0)$ is an algebraic integer for all the modular functions $F(z)=\frac{\Theta_M(D'^e\frac{z}{2^{e}3^{e'}m^{e''}})}{\Theta_K(\frac{z}{m^{e'''}})}$ defined above in Lemma \ref{shim_1} and Lemma \ref{shim_3}. This is a standard result from CM theory: for $\gamma \in \SL_2(\ZZ)$, each $F(\gamma z)$ is a modular function with algebraic coefficients in its Fourier expansion; then $F(\tau_0)$ is an algebraic integer when evaluated at a CM point $\tau_0$. Thus $F(\tau_0)D^{1/3}\alpha^{1/2}$ are algebraic integers and their traces are in $\OO_K$. Thus we can state:
\begin{Prop}\label{OK}
\begin{enumerate}[(i)]
	\item $X^*, U, T, Y_i, Y, W_i\in \OO_K$.
	\item $Z, Y_i^{\circ}, Y^{\circ}, W_i^{\circ}, U^{\circ} \in \OO_K$
\end{enumerate}
\end{Prop}

\bigskip
\subsection{Important lemma} The following Lemma is essential in computing the complex conjugates of $Y_0$ and $W_1$:

\begin{Lem}\label{lemma_A} Let $\tau=\frac{-b_0+\sqrt{-3}}{2}$ such that $b_0 \equiv \sqrt{-3} (\alpha)$. For the ideal $\A_0=(k_0)$ ideal generated by $k_0=t_{\A}a+s_{\A}\tau$ of norm $a$ such that $3m^*|s_{\A}$ and $D'|t_{\A}$, we have: 
\[
\frac{\Theta_M(D'\frac{\tau}{3^{e} m^{e'}})}{\Theta_K(\frac{\tau}{3^{e''}})}
=
\frac{1}{D'}\left(\frac{D}{3m^*}\right)\left(\frac{t_{\A}}{m^*}\right)\left(\frac{\Theta_M(\frac{\tau}{D'3^{e} m^{e'}})}{\Theta_K(\frac{\tau}{3^{e''}})}\right)^{\sigma_{\A_{0}}^{-1}},
\]
where $e, e', e''\in\{0, 1\}$.
\end{Lem}

{\bf Proof:} Let $c=\frac{b_0^2+3}{4am}$. Note that $ G(z)=\frac{\Theta_M(\frac{z}{3^{e} m^{e'}})}{\Theta_K(z/3^{e''})}$ is a modular function of level dividing $6D'm^*$ and, from Shimura reciprocity law \eqref{Shimura}, we have:
\[
G(\tau)^{\sigma_{\A_0}^{-1}}
=
G(\left(\begin{smallmatrix} t_{\A}a-s_{\A}b & -s_{\A}cm \\ s_{\A} & t\end{smallmatrix}\right)\tau).
\]
Explicitly, we compute:
\begin{equation}\label{comp1}
\Theta_M(\left(\begin{smallmatrix} D' & 0 \\ 0 & 3^{e}m^{e'}\end{smallmatrix}\right)\left(\begin{smallmatrix} t_{\A}a-s_{\A}b & -s_{\A}cm \\ s & t\end{smallmatrix}\right)\tau)
=
\Theta_M(\left(\begin{smallmatrix} (t_{\A}a-s_{\A}b)D' & -s_{\A}cm^{1-e'}/3^{e} \\ -s_{\A}m^{e'}3^{e} & t_{\A}/D'\end{smallmatrix}\right)\left(\begin{smallmatrix} 1 & 0 \\ 0 & D'3^em^{e'}\end{smallmatrix}\right)\tau).
\end{equation}

As $s_{\A}=m^*s$ with $3|s$, we have $3m^*|-s_{\A}m^{e'}3^{e}$, and using Lemma \ref{modular}, we get in \eqref{comp1} $\left(\frac{t_{\A}/D'}{3m^*}\right)\frac{1}{D'}(s_{\A}\tau+t_{\A})\Theta_M(\frac{\tau}{D'3^em^{e'}})$.
 Similarly we get
 \[
 \Theta_K(\left(\begin{smallmatrix} 1 & 0 \\ 0 & 3^{e''}\end{smallmatrix}\right)\left(\begin{smallmatrix} t_{\A}a-s_{\A}b & -s_{\A}cm \\ s_{\A} & t_{\A}\end{smallmatrix}\right)\tau)
 =
\left(\frac{t_{\A}}{3}\right) (s_{\A}\tau+t_{\A})\Theta_K(\tau).
 \]

Thus we have $\ds G(\left(\begin{smallmatrix} t_{\A}a-s_{\A}b & -s_{\A}cm \\ s_{\A} & t\end{smallmatrix}\right)\tau_0)
=
\frac{1}{D'}\left(\frac{D}{3m^*}\right)\left(\frac{t_{\A}}{m^*}\right)\frac{\Theta_M(\frac{\tau}{D'3^em^{e'}})}{\Theta_K(\tau/3^{e''})}.$

\bigskip

\subsection{Condition for $m\equiv 3(4)$}\label{XZ_i}

We recall $b_i\equiv \sqrt{-3}(\alpha)$ such that $b_1\equiv 1(8)$ and $b_3\equiv -1(8)$, and $\tau_i=\frac{-b_i+\sqrt{-3}}{2}$. We define for $m\equiv 3(4)$:
\begin{itemize}
	\item $X^{(i)}=\Tr_{H_{6Dm}/K} \frac{\Theta_M\left(D\tau_i/2\right)}{\Theta_K (\omega)} \alpha^{1/2}D^{1/3}$, $X=\Tr_{H_{12Dm}/K} \frac{\Theta_M\left(D\tau_i\right)}{\Theta_K (\omega)} \alpha^{1/2}D^{1/3}$
	
	\item $Z'^{(i)}=\Tr_{H_{6D}/K} \frac{\Theta_M\left(D\tau_i/(2m)\right)}{\Theta_K (\tau_i/m)} \alpha^{1/2}D^{1/3}$, $Z'=\Tr_{H_{12D}/K} \frac{\Theta_M\left(D\tau_i/m\right)}{\Theta_K (\tau_i/m)} \alpha^{1/2}D^{1/3}$
	
\end{itemize}

We note that $Z'$ and $Z'^{(i)}$ are well defined from Lemma \ref{shim_3} and they take values in $\OO_K$.

\begin{Lem}\label{X_zero} For $\alpha \equiv \pm \sqrt{-3} (4)$, we have $X^{(i)}=0$ and $Z'^{(i)}=0$ for $\left[\alpha\right]\left(\frac{b_i}{4}\right)=1$. 
\end{Lem}

{\bf Proof:} Let $\tau_i=\frac{-b_i+\sqrt{-3}}{2}$ and let $\A$ be an ideal prime to $6mD$, with generator $\beta=ta+sm\tau_i$ such that $2||s$, $3D|s$. 
We note that $\Nm\A=a\equiv 3(4)$ and we fix the generator of $\A$ such that $\beta\equiv \sqrt{-3} (4)$. Note that then $t\equiv -b_i(4)$.

Then we have by Shimura reciprocity law \eqref{Shimura_A}, for the modular function $F(z)=\Theta(Dz/2)/\Theta(z)$ of level $6Dm$, $F(\tau_0)=F(\left(\begin{smallmatrix}ta-sb & -sc \\ s & t \end{smallmatrix}\right)\tau_i).$ Explicitly, we compute
\[
\Theta_M(\left(\begin{smallmatrix}D & 0 \\ 0 & 2 \end{smallmatrix}\right)\left(\begin{smallmatrix}ta-sb & -smc \\ sm & t \end{smallmatrix}\right)\tau_i)=\Theta_M(\left(\begin{smallmatrix}ta-sb & -Dsmc/2 \\ 2sm/D & t \end{smallmatrix}\right)(D\tau_i/2)).
\]
We note that $12m| 2sm/D$, thus the matrix above is in $\Gamma_0(3m^*)$ and we can apply Lemma \ref{modular}. Thus we get above $(sm\tau_i+t)\left(\frac{t}{3m^*}\right)\Theta_M(D\tau_i/2)$. Similarly we compute $\Theta_K(\left(\begin{smallmatrix}ta-sb & -smc \\ sm & t \end{smallmatrix}\right)\tau_i)=(sm\tau_i+t)\left(\frac{t}{3}\right)\Theta_K(\tau_i)$. Thus we get:
\[
F(\tau_i)^{\sigma_{\A}^{-1}}=\left(\frac{t}{m^*}\right)F(\tau_i).
\]
On the other hand $(\alpha^{1/2})^{\sigma_{\A}^{-1}}=\alpha^{1/2}\left[\frac{\alpha}{\A}\right]$. From the reciprocity law \ref{rec_law0}, we have $\left[\frac{\alpha}{\A}\right]=\left[\frac{ta}{\alpha}\right][\overline{\alpha}][\beta][\overline{\alpha}\beta]$. We compute $\left[\frac{ta}{\alpha}\right]=\left(\frac{ta}{m}\right)$ and $[\overline{\alpha}][\beta][\overline{\alpha}\beta]=-[\overline{\alpha}]=[\alpha]$, as $\beta\equiv \sqrt{-3}(4)$. Thus we get 
\[
(F(\tau_i)\alpha^{1/2}D^{1/3})^{\sigma_{\A}^{-1}}=\left(\frac{t}{4}\right)[\alpha]F(\tau_i)\alpha^{1/2}D^{1/3}
\]
which equals $-\left(\frac{b_i}{4}\right)[\alpha]F(\tau_i)\alpha^{1/2}D^{1/3}$. Thus, if $\left(\frac{b_i}{4}\right)[\alpha]=1$, the trace of $X^{(i)}$ equals $0$.

The proof is similar for $Z'^{(i)}$, by taking $\A'$ be an ideal prime to $6D$, with generator $\beta=ta+s\tau_i$ such that $2||s$, $3D|s$. Then $c=\frac{b_i^2+3}{4}$ is divisible by $m$.  For the modular function $F'(z)=\Theta(\frac{Dz}{2m})/\Theta(\frac{z}{m})$ of level $6D$, we have $F'(\tau_0)=F'(\left(\begin{smallmatrix}ta-sb & -scm \\ s & t \end{smallmatrix}\right)\tau_i)$ which gives us:
\[
\Theta_M(\left(\begin{smallmatrix}D & 0 \\ 0 & 2m \end{smallmatrix}\right)\left(\begin{smallmatrix}ta-sb & -scm \\ s & t \end{smallmatrix}\right)\tau_i)=\Theta_M(\left(\begin{smallmatrix}ta-sb & -Dsc/2 \\ 2sm/D & t \end{smallmatrix}\right)\frac{D\tau_i}{2m}),
\]
which equals $\left(\frac{t}{3m^*}\right)(s\tau_i+t)\Theta_M(D\tau_i/(2m))$. Similarly we get $\Theta_K(\left(\begin{smallmatrix}1& 0 \\ 0 & m \end{smallmatrix}\right)\left(\begin{smallmatrix}ta-sb & -scm \\ s & t \end{smallmatrix}\right)\tau_i)=\left(\frac{t}{3}\right)(s\tau_i+t)\Theta_K(\tau_i/m)$, thus $F'(\tau_i)^{\sigma_{\A'}^{-1}}=\left(\frac{t}{m^*}\right)F(\tau_i)$, and again $(F'(\tau_i)\alpha^{1/2}D^{1/3})^{\sigma_{\A'}^{-1}}=-\left(\frac{b_i}{4}\right)[\alpha]F'(\tau_i)\alpha^{1/2}D^{1/3}.$ This implies $Z'^{(i)}=0$ under the condition $\left(\frac{b_i}{4}\right)[\alpha]=1$.

\begin{cor}\label{X_1} For $\alpha \equiv \sqrt{-3} (4)$, for $i$ such that $\left[\alpha\right]\left(\frac{-b_i}{4}\right)=1$, we have $X=X^{(i)}$ and $Z'=Z'^{(i)}$.

\end{cor}

{\bf Proof:} For $\tau_i=\frac{-b_i+\sqrt{-3}}{2}$ with $b_i\equiv i(4)$,  from Lemma \ref{sum_M} for $d=2$, we have $\Theta_M(D\tau_1/2)+\Theta_M(D\tau_3/2)=2\Theta_M(D\tau_i)$, thus by multiplying by $\alpha^{1/2}D^{1/3}/\Theta_K(\omega)$ and taking the traces from $H_{12mD}$ to $K$, we get:
\[
2X^{(1)}+2X^{(3)}=2X.
\]
However, since $X^{(3)}=0$ for $\alpha\equiv \sqrt{-3}$ and $X^{(1)}=0$ for $\alpha\equiv -\sqrt{-3}$ from Lemma \ref{X_zero}, we get the result. The proof is similar for $Z'$.

\begin{rmk} We note that for $b_i$ under the condition \eqref{cond}, we have $Z'^{(i)}=U^{\circ}$, thus $Z'=Z'^{(i)}=U^{\circ}$.

\end{rmk}

\bigskip

\section{Appendix}

\subsection{Tamagawa numbers}\label{tam} Using Tate's algorithm (see \cite{S2}), one can compute the Tamagawa number $c_{E_{D, \alpha}}=c_2c_3c_Dc_{\alpha} $, where:
	\begin{itemize}
		\item $c_2=1$, 						\item $c_3=\begin{cases} 3 & D\equiv \pm 1(9), \left[\frac{\alpha}{\sqrt{-3}}\right]=1$ or $D\equiv \pm 4(9), \left[\frac{\alpha}{\sqrt{-3}}\right]=-1\\
1 & D\equiv \pm 1(9), \left[\frac{\alpha}{\sqrt{-3}}\right]=-1$ or  $D\equiv \pm 4(9), \left[\frac{\alpha}{\sqrt{-3}}\right]=1\\
4& D\equiv \pm 2(9)\end{cases}$,
	\item $c_D=3^{\#\{\pp|D: \left(\frac{\alpha}{\pp}\right)=1\}}$, which gives us $c_{D}=\begin{cases} 3^{2k} & \left(\frac{D}{3m^*}\right)=1\\
 3^{2k+1} & \left(\frac{D}{3m^*}\right)=-1\end{cases}$.

	\item $c_{\alpha}=\begin{cases}1 & 2D^2\not\equiv u^3(\alpha) \\ 4 & -2D^2\equiv u^3(\alpha) \end{cases}=\begin{cases}1 & \chi_{2D^2}(\alpha)\neq 1 \\ 4 & \chi_{2D^2}(\alpha)=1 \end{cases}$.
  \end{itemize}
Thus $c_0=(-1)^{v_3(\prod_{p|6Dm}c_p)}=\begin{cases} -\left(\frac{D}{3m^*}\right)\left[\frac{\alpha}{\sqrt{-3}}\right], & D\equiv \pm 1\\ \left(\frac{D}{3m^*}\right)\left[\frac{\alpha}{\sqrt{-3}}\right], & D\equiv \pm 4 (9) \\ \left(\frac{D}{3m^*}\right), & D\equiv \pm 2(9) \end{cases}.$

\subsection{Computations characters}

We recall the constant $u_{\alpha, b}=\begin{cases} 1&   \text{ if } \alpha\equiv 1(4)\\-1&   \text{ if } \alpha\equiv -1(4) \\ [\alpha]\left(\frac{-b}{4}\right)&  \text{ if } \alpha\equiv \pm\sqrt{-3}(4) \end{cases}$.	

We also recall the ideals 
\begin{itemize}
	\item $\C_0=\left(\frac{-b_{0, 0}+\sqrt{-3}}{2}\right)/(\alpha\sqrt{-3})$ of norm $c_0$, with $b_{0, 0} \equiv 0(3)$, $b_{0, 0}\equiv\sqrt{-3}(\alpha)$, $b_{0, 0}\equiv -1(2D)$
	\item $\C_i=\left(\frac{-b_{0, i}+\sqrt{-3}}{2}\right)/(\alpha)$ of norm $c_i$, with $b_{0, i} \equiv -i(3)$, $b_{0, i}\equiv\sqrt{-3}(\alpha)$, $b_{0, i}\equiv -1(2D)$, for $i=1, 2$
	\item $\A_0=t_{\A}a+3m^*s'\frac{-b_{0, 0}+\sqrt{-3}}{2}$, with $D' |t_{\A}$, of norm $a_0$
	\item $\A_j=\left(\frac{-b_j+\sqrt{-3}}{2}\right)$ of norm $a_j$ with $b_j-b_0\equiv j (m)$, $b_j^2\neq -3(m)$ and $b_j\equiv 1(6D)$ 
	\item $\C_0^{*}=(\frac{-b_0^*+\sqrt{-3}}{2})/(\overline{\alpha})$ with $b_0^*\equiv -\sqrt{-3}(\alpha)$, $b_0^*\equiv 1(6D)$
\end{itemize}
We also assume that all $b_{0, i}, b_0^*, b_j$ are $\equiv \pm 1(8)$. We compute below the character $\eps$ for each of the ideals defined:

\begin{Lem}\label{char} \begin{enumerate}[(i)]
	\item $\eps(\overline{\C_0})=u_{\alpha, b_{0, 0}}\eps(\overline{\alpha})$
	
	\item $\eps(\overline{\C_i})=u_{\alpha, b_{0, i}}\eps(\overline{\alpha})\eps(\sqrt{-3})$
	
	\item $\eps(\overline{\A_0})=\left(\frac{t}{m^*}\right)$
	
	\item $\eps(\A_j)=\begin{cases} \left(\frac{\left((-b_j+b_0)/2\right)}{m}\right)\left(\frac{b_j}{4}\right) &  \text{for } \alpha\equiv \pm \sqrt{-3}(4)  \\ \left(\frac{\left((-b_j+b_0)/2\right)}{m}\right)&  \text{for } \alpha\equiv \pm 1(4) \end{cases}$
		
	\item $\eps(\C^*_0)=\eps(\overline{\alpha})\eps(\sqrt{-3})u_{\alpha, -b_0^*}$

\end{enumerate}

\end{Lem}

{\bf Proof:} For $(i), (ii)$ and $(v)$, we first note that for $b\equiv \sqrt{-3}(\alpha)$: $\eps(\frac{b+\sqrt{-3}}{2})=\left[\frac{\alpha}{\frac{b+\sqrt{-3}}{2}}\right]=\left[\frac{\frac{b+\sqrt{-3}}{2}}{\alpha}\right][\alpha]\left[(\frac{b-\sqrt{-3}}{2})\right]\left[\alpha(\frac{b-\sqrt{-3}}{2})\right]$
 from the reciprocity law \eqref{rec_law0}. This equals $\left[\frac{\sqrt{-3}}{\alpha}\right][\alpha]\left[(\frac{b-\sqrt{-3}}{2})\right]\left[\alpha(\frac{b-\sqrt{-3}}{2})\right]$. We further apply the reciprocity law $\left[\frac{\sqrt{-3}}{\alpha}\right]=\left[\frac{\alpha}{\sqrt{-3}}\right][\alpha][-\sqrt{-3}][-\alpha\sqrt{-3}]=\eps(\sqrt{-3})[\alpha][-\alpha\sqrt{-3}]$. Thus $\eps(\frac{b+\sqrt{-3}}{2})=\eps(\sqrt{-3})[\alpha]\left[(\frac{b-\sqrt{-3}}{2})\right]\left[\alpha(\frac{b-\sqrt{-3}}{2})\right][-\alpha\sqrt{-3}]$, which equals:
 \[
 \eps\left(\frac{b+\sqrt{-3}}{2}\right)= \eps(\sqrt{-3})u_{\alpha, b}
 \] for $b\equiv \pm 1(8)$. Then we can compute:

\begin{enumerate}[(i)]
	\item $\eps(\overline{\alpha})\eps(\overline{\C_0})\eps(\sqrt{-3})=\eps(\frac{b_{0, 0}+\sqrt{-3}}{2})$, thus $\eps(\overline{\C_0})=\eps(\overline{\alpha})u_{\alpha, b_{0, 0}}$.

	\item $\eps(\overline{\C_i})=\eps(\overline{\alpha})\eps(\frac{b_{0, i}+\sqrt{-3}}{2})$, thus $\eps(\overline{\C_i})=\eps(\overline{\alpha})\eps(\sqrt{-3})u_{\alpha, b_{0, i}}$.

	\item As $\eps(\overline{\A_0})=\left[\frac{\alpha}{\overline{A_0}}\right]$, the reciprocity law \eqref{rec_law0} gives us $\left[\frac{\alpha}{\A_0}\right]=\left[\frac{ta}{\alpha}\right][\alpha][ta][ta\alpha]$. This equals $\left(\frac{ta}{m}\right)$ for $m\equiv 1(4)$. As $a\equiv1/t^2_{\A}(m)$, we get $\left(\frac{t}{m}\right)$.
	
	For   $\alpha\equiv \pm \sqrt{-3}(4)$, we get $\left(\frac{ta}{m}\right)\left(\frac{ta}{4}\right)$ and, similarly,  $a\equiv1/t^2_{\A}(m^*)$, thus we get $\left(\frac{t}{4m}\right)$.

	\item  $\eps(\A_i)=\left[\frac{\alpha}{(-b_i/2+\sqrt{-3}/2)}\right]$. We get $\left[\frac{\alpha}{\left(-b_i/2+\sqrt{-3}/2\right)}\right]=\left[\frac{\left(-b_i/2+b_0/2\right)}{\alpha}\right][\overline{\alpha}]\left[\left(\frac{-b_i+\sqrt{-3}}{2}\right)\right][\overline{\alpha}\left(\frac{-b_i+\sqrt{-3}}{2}\right)]$. This equals $\left(\frac{(-b_i+b_0)/2}{m}\right)\left(\frac{b_i}{4}\right)$ for $\alpha \equiv \pm \sqrt{-3} (4)$ and $\left(\frac{(-b_i+b_0)/2}{m}\right)$ for $\alpha \equiv \pm 1(4)$.

	\item $\eps(\C^*_0)\eps(\overline{\alpha})=\eps(\frac{-b_0^*+\sqrt{-3}}{2})$, thus $\eps(\C^*_0)=\eps(\overline{\alpha})\eps(\sqrt{-3})u_{\alpha, -b_0^*}$.

	\end{enumerate}

\bigskip

We also compute the $\chi_D$ values: 

\begin{Lem}\label{charD} We have
\begin{enumerate}[(i)]
	\item $\chi_D(\overline{\C_0})=\chi_D(\left(\frac{b_{0, 0}/3\sqrt{-3}+1}{2}\right))\chi_D(\alpha)$
	
	\item $\chi_D(\overline{\C_1})=\chi_D(\overline{\alpha})$,
	$\chi_D(\overline{\C_2})=\chi_{\omega}(D)\chi_D(\overline{\alpha})$
	
	\item $\chi_D(\overline{\A_0})=\chi_{\omega}(D)$
	
	\item $\chi_D(\overline{\A_i})=1$

	\item $\chi_D(\overline{\C_0}^*)=1$

\end{enumerate}
\end{Lem}

{\bf Proof:} \begin{enumerate}[(i)]
	
	\item $\chi_D(\C_0)\chi_D(\alpha)=\chi_D(\left(\frac{b_{0, 0}/3\sqrt{-3}+1}{2}\right))$.

	 \item $\chi_D(\overline{\C_i})=\chi_D(\alpha)\chi_D( \frac{b_{0, i}+\sqrt{-3}}{2})$. For $i=1$ we have $\left(\frac{b_{0, i}+\sqrt{-3}}{2}\right)\equiv \omega (3D)$, thus $\chi_D(\left(\frac{b_{0, 1}+\sqrt{-3}}{2}\right))=1$. For $i=2$, $\left(\frac{b_{0, 2}+\sqrt{-3}}{2}\right)\equiv \omega (D)$, $\equiv -\omega^2 (3)$, thus $\chi_D(\left(\frac{b_{0, 2}+\sqrt{-3}}{2}\right))=\chi_{\omega}(D)$

	\item As $t_{\A}a+3s'm\frac{-b_{0}-\sqrt{-3}}{2}\equiv \pm1(3)$ and $D|t_{\A}$, $\frac{-b_{0}-\sqrt{-3}}{2}\equiv -\omega(D)$, thus $\chi_D(\overline{\A_0})=\chi_{\omega}(D)$
	
	\item $\frac{-b_i+\sqrt{-3}}{2}\equiv \omega (3D)$, thus $\chi_D(\overline{\A_i})=1$.
	
	\item $\tau_0^*\equiv \omega(3D)$, thus $\chi_D(\C_0^*)=1$.
	\end{enumerate}

\subsection{Computing Galois conjugates}

The following Lemma is used in Sections \ref{XT} and \ref{XYZ} to show that various terms are Galois conjugate to each other. The proof consists of applying Shimura reciprocity law \eqref{Shimura_A} and the inverse transformations \eqref{transform} and \eqref{transform2} for the theta functions $\Theta_M$ and $\Theta_K$, adjusting the formulas using \eqref{lemma_phi}.
 
\begin{Lem}\label{inv1} For the Galois actions corresponding to the ideals $\C_i, \C_0^*, \C_0', \A_i$ defined above, we have:
\begin{enumerate}[(i)]
	\item $\frac{\Theta_M\left(-\overline{\tau_{0,0}}/3D'm\right)}{\Theta_K\left(\omega\right)}
	=
	\frac{3\varphi(\alpha) D'}{\sqrt{m}}\left(\frac{\Theta_M\left(D'\tau_{0, 0}/3m\right)}{\Theta_K\left(\tau_0/9\right)}\right)^{\sigma_{\C_0}^{-1}}$ for $m\equiv 1(4)$

$\frac{\Theta_M\left(-\overline{\tau_{0,0}}/6Dm\right)}{\Theta_K\left(\omega\right)}=\frac{3\varphi(\alpha) D}{\sqrt{m}}\left(\frac{\Theta_M\left(D\tau_{0,0}/6m\right)}{\Theta_K\left(\tau_0/9\right)}\right)^{\sigma_{\C_0}^{-1}}$ for $m\equiv 3(4)$
	
	\item $
\frac{\Theta_M\left(-\overline{\tau_{0, i}}/3mD'\right)}{\Theta_K\left(-\overline{\tau_{0, i}}/3\right)}
=\frac{D'\sqrt{m}}{\varphi(\overline{\alpha})}\left(\frac{\Theta_M\left(D'\tau_{0, i}/m\right)}{\Theta_K(\omega)}\right)^{\sigma_{\C_i}^{-1}}$ for $m\equiv 1(4)$

 $\frac{\Theta_M\left(-\overline{\tau_{0, i}}/6mD\right)}{\Theta_K\left(-\overline{\tau_{0, i}}/3\right)}
=\frac{D\sqrt{m}}{\varphi(\overline{\alpha})}\left(\frac{\Theta_M\left(D\tau_{0, i}/2m\right)}{\Theta_K(\omega)}\right)^{\sigma_{\C_i}^{-1}}$ for $m\equiv 3(4)$

\item $\frac{\Theta_M(-\overline{\tau_i}/3mD')}{\Theta_M(-\overline{\tau_i}/3)}
	=
	D'\sqrt{m}\left(\frac{\Theta_M(D'\omega)}{\Theta_M(\omega)}\right)^{\sigma_{\A_i}^{-1}}$ for $m\equiv 1(4)$
	
 $\frac{\Theta_M(-\overline{\tau_i}/6mD)}{\Theta_M(-\overline{\tau_i}/3)}	=
	D\sqrt{m}\left(\frac{\Theta_M(D\tau_i/2)}{\Theta_M(\omega)}\right)^{\sigma_{\A_i}^{-1}}$ for $m\equiv 3(4)$

	\item $\frac{\Theta_M(-\overline{\tau_0}/3mD')}{\Theta_M(-\overline{\tau_0}/3)}
=
\frac{D'\sqrt{m}}{\varphi(\overline{\alpha})}\left(\frac{\Theta_M(D'\tau_0/m)}{\Theta_M(\omega)}\right)^{\sigma_{\C'_0}^{-1}} 
$  for $m\equiv 1(4)$

 $\frac{\Theta_M(-\overline{\tau_0}/6mD)}{\Theta_M(-\overline{\tau_0}/3)}
=
\frac{D\sqrt{m}}{\varphi(\overline{\alpha})}\left(\frac{\Theta_M(D\tau_0/2m)}{\Theta_M(\omega)}\right)^{\sigma_{\C'_0}^{-1}}$ for $m\equiv 3(4)$

	\item $\frac{\Theta_M(-\overline{\tau_0}^*/3mD')}{\Theta_M(-\overline{\tau_0}^*/3)}
	=
\frac{D'\sqrt{m}}{\varphi(\alpha)}\left(\frac{\Theta_M(D'\tau_0^*/m)}{\Theta_M(\omega)}\right)^{\sigma_{\C^*_0}^{-1}}$ for $m\equiv 1(4)$ 

$\frac{\Theta_M(-\overline{\tau_0}^*/6mD')}{\Theta_M(-\overline{\tau_0}^*/3)}
	=
\frac{D\sqrt{m}}{\varphi(\alpha)}\left(\frac{\Theta_M(D\tau_0^*/2m)}{\Theta_M(\omega)}\right)^{\sigma_{\C^*_0}^{-1}}$ for $m\equiv 3(4)$

\item $\frac{\Theta_M(D'\tau_i/m)}{\Theta_K(\omega)}=
	\frac{\sqrt{m}}{D'}\left(\frac{\Theta_M(-\overline{\tau_i}/3D')}{\Theta_K(-\overline{\tau_i}/3)}\right)^{\sigma_{\overline{\A_i}}^{-1}}$ for $m\equiv 1(4)$
	
$\frac{\Theta_M(D\tau_i/2m)}{\Theta_K(\omega)}=
	\frac{\sqrt{m}}{D}\left(\frac{\Theta_M(-\overline{\tau_i}/6D')}{\Theta_K(-\overline{\tau_i}/3)}\right)^{\sigma_{\overline{\A_i}}^{-1}}$ for $m\equiv 3(4)$

\end{enumerate}

\end{Lem}

{\bf Proof:} We sketch below the proofs for $m\equiv 1(4)$. 
\begin{enumerate}[(i)]
	\item $\frac{\Theta_M\left(\frac{b_{0,0}+\sqrt{-3}}{6D'm}\right)}{\Theta_K\left(\frac{b_{0, 0}+\sqrt{-3}}{2}\right)}
=
\varphi(\alpha) \frac{\Theta_M\left(\frac{b_{0,0}+\sqrt{-3}}{6D'm}\right)}{\Theta_K\left(\frac{b_{0,0}+\sqrt{-3}}{2m}\right)}
=
\frac{3\varphi(\alpha)D'}{\sqrt{m}} \frac{\Theta_M\left(\frac{-b_{0, 0}+\sqrt{-3}}{6c_0m}D'\right)}{\Theta_K\left(\frac{-b_{0, 0}+\sqrt{-3}}{18c_0}\right)}
=
\frac{3\varphi(\alpha) D'}{\sqrt{m}}\left(\frac{\Theta_M\left(D'\frac{-b_{0, 0}+\sqrt{-3}}{6m}\right)}{\Theta_K\left(\frac{-b_{0, 0}+\sqrt{-3}}{18}\right)}\right)^{\sigma_{\C_0}^{-1}},$
where we first have used \eqref{lemma_phi} for the ideal $\alpha=[m, \frac{-b_{0, 0}+\sqrt{-3}}{2}]$, followed by the transformations \eqref{transform} and \eqref{transform2}, and finally using Shimura reciprocity \eqref{Shimura_A} for the ideal $\C_0=\left(\frac{-b_{0, 0}+\sqrt{-3}}{2}\right)/(\alpha\sqrt{-3})$ of norm $c_0$.

\item $\frac{\Theta_M\left(\frac{b_{0, i}+\sqrt{-3}}{6mD'}\right)}{\Theta_K\left(\frac{b_{0, i}+\sqrt{-3}}{6}\right)}
= D'\sqrt{m}\frac{\Theta_M\left(D'\frac{-b_{0, i}+\sqrt{-3}}{2mc_i}\right)}{\Theta_K\left(\frac{-b_{0, i}+\sqrt{-3}}{2mc_i}\right)}=D'\sqrt{m}
\left(\frac{\Theta_M\left(D'\frac{-b_{0, i}+\sqrt{-3}}{2m}\right)}{\Theta_K\left(\frac{-b_{0, i}+\sqrt{-3}}{2m}\right)}\right)^{\sigma_{\C_i}^{-1}}
=
\frac{D'\sqrt{m}}{\varphi(\overline{\alpha})}\left(\frac{\Theta_M\left(D'\frac{-b_{0, i}+\sqrt{-3}}{2m}\right)}{\Theta_K(\omega)}\right)^{\sigma_{\C_i}^{-1}},$
where we have used \eqref{transform} and \eqref{transform2}, followed by Shimura reciprocity  \eqref{Shimura_A} for the ideal
 $\C_i=(\frac{-b_{0, i}+\sqrt{-3}}{2})/(\alpha)$ of norm $c_i$, and in the last step we applied \eqref{lemma_phi}.

\item $\frac{\Theta_M(\frac{b_i+\sqrt{-3}}{6mD'})}{\Theta_M(\frac{b_i+\sqrt{-3}}{6})}
	=
	D'\sqrt{m}\left(\frac{\Theta_M(D'\frac{-b_i+\sqrt{-3}}{2a_i})}{\Theta_M(\frac{-b_i+\sqrt{-3}}{2a_i})}\right)
	=
	D'\sqrt{m}\left(\frac{\Theta_M(D'\frac{-b_i+\sqrt{-3}}{2})}{\Theta_M(\frac{-b_i+\sqrt{-3}}{2})}\right)^{\sigma_{\A_i}^{-1}}$, where we have used \eqref{transform} and \eqref{transform2}, followed by \eqref{Shimura_A} for the ideal $\A_i=\left(\frac{-b_i+\sqrt{-3}}{2}\right)$ with norm $a_i=\frac{b_i^2+3}{2}$.

	\item $\frac{\Theta_M(\frac{b_0+\sqrt{-3}}{6mD'})}{\Theta_M(\frac{b_0+\sqrt{-3}}{6})}
=
\frac{D'\sqrt{m}}{\varphi(\overline{\alpha})}\left(\frac{\Theta_M(D'\frac{-b_0+\sqrt{-3}}{2m})}{\Theta_M(\frac{-b_0+\sqrt{-3}}{2})}\right)^{\sigma_{\A_0}^{-1}}
$ 
where we first used \eqref{transform} and \eqref{transform2}, followed by \eqref{Shimura_A} for $\A_0=\left[a_0, \frac{-b_0+\sqrt{-3}}{2}\right]=(\frac{-b_0+\sqrt{-3}}{2})/(\alpha)$, and finally \eqref{lemma_phi}.

\item $\frac{\Theta_M(\frac{b_0^*+\sqrt{-3}}{6mD'})}{\Theta_M(\frac{b^*_0+\sqrt{-3}}{6})}
	=
	D'\sqrt{m}\frac{\Theta_M(D'\frac{-b_0^*+\sqrt{-3}}{2mc_0^*})}{\Theta_M(\frac{-b_0^*+\sqrt{-3}}{2mc_0^*})}
	=
D'\sqrt{m}\left(\frac{\Theta_M(D'\frac{-b_0^*+\sqrt{-3}}{2m})}{\Theta_M(\frac{-b_0^*+\sqrt{-3}}{2m})}\right)^{\sigma_{\C^*_0}^{-1}}
=
\frac{D'\sqrt{m}}{\varphi(\alpha)}\left(\frac{\Theta_M(D'\frac{-b_0+\sqrt{-3}}{2m})}{\Theta_M(\frac{-b_0+\sqrt{-3}}{2})}\right)^{\sigma_{\C^*_0}^{-1}}$, where we have used \eqref{transform} and \eqref{transform2}, followed by Shimura reciprocity  \eqref{Shimura_A} for the ideal
 $\C_0^*=(\frac{-b_{0}^*+\sqrt{-3}}{2})/(\overline{\alpha})$ of norm $c_0^*$, and in the last step we applied \eqref{lemma_phi}.

\item $\frac{\Theta_M(D'\frac{-b_i+\sqrt{-3}}{2m})}{\Theta_K(\frac{-b_i+\sqrt{-3}}{2})}=\frac{\sqrt{m}}{D}\frac{\Theta_M(\frac{b_i+\sqrt{-3}}{6D'a_i})}{\Theta_K(\frac{b_i+\sqrt{-3}}{6a_i})}
	=
	\frac{\sqrt{m}}{D'}\left(\frac{\Theta_M(\frac{b_i+\sqrt{-3}}{6D'})}{\Theta_K(\frac{b_i+\sqrt{-3}}{6})}\right)^{\sigma_{\overline{\A_i}}^{-1}}$, where we have first used \eqref{transform} and \eqref{transform2}, followed by \eqref{Shimura_A} for the ideal $\overline{\A_i}=\left(\frac{b_i+\sqrt{-3}}{2}\right)$ of norm $a_i$.

\end{enumerate}
The proofs for $m\equiv 3(4)$ are completely similar.

\begin{small}

\end{small}

\end{document}